\documentclass[12pt]{amsart}
\usepackage{amssymb}
\usepackage{amsfonts}
\usepackage{latexsym}
\usepackage{amscd}

\newcommand{\N}{\mathbb{N}}
\newcommand{\Z}{\mathbb{Z}}

\newcommand{\C}{\mathbb{C}}

\newtheorem{lemma}{Lemma}[section]
\newtheorem{corollary}[lemma]{Corollary}
\newtheorem{theorem}[lemma]{Theorem}
\newtheorem{proposition}[lemma]{Proposition}

\newtheorem{definition}[lemma]{Definition}
\newtheorem{example}[lemma]{Example}

\usepackage[all]{xy}

\input xy
\xyoption{all}

\begin{document}
\title[Isomorphisms between Leavitt algebras]{Isomorphisms between Leavitt algebras and their matrix rings}
\author{G. Abrams}
\address{Department of Mathematics, University of Colorado,
Colorado Springs CO 80933 U.S.A.} \email{abrams@math.uccs.edu}
\author{P. N. \'{A}nh}
\address{R\'enyi Institute of Mathematics, Hungarian Academy of
Sciences, 1364 Budapest, Pf. 127 Hungary} \email{anh@renyi.hu}
\author{E. Pardo}
\address{Departamento de Matem\'aticas, Universidad de C\'adiz,
Apartado 40, 11510 Puerto Real (C\'adiz), Spain.}\email{enrique.pardo@uca.es}
\urladdr{http://www2.uca.es/dept/matematicas/PPersonales/PardoEspino/index.HTML}
\thanks{The first author is grateful for support provided by the Clear Creek Ranch Institute and its directors,
L. George and S. Hyde.}
 \thanks{The second author is supported partly by Hungarian National Foundation for Scientific Research grant no. K61007.
  During Fall 2006 he was also supported partly by The Colorado College, the University of Colorado at Colorado Springs, and Professor Michael
  Siddoway.}
\thanks{The
third author was partially supported by the DGI and European Regional Development Fund, jointly, through Project
MTM2004-00149, by PAI III projects FQM-298 and P06-FQM-1889 of the Junta de Andaluc\'{\i}a, by the Comissionat per
Universitats i Recerca de la Generalitat de Catalunya, and by the Consolider Ingenio
``Mathematica" project CSD2006-32 by the MEC} 
\date{18 December 2006}
%
%


\begin{abstract}
Let $K$ be any field, let $L_n$ denote the Leavitt algebra of type $(1,n-1)$ having coefficients in $K$, and let ${\rm
M}_d(L_n)$ denote the ring of $d \times d$ matrices over $L_n$. In our main result, we show that ${\rm M}_d(L_n) \cong
L_n$ if and only if $d$ and $n-1$ are coprime. We use this isomorphism to answer a question posed in \cite{PS}
regarding isomorphisms between various C*-algebras.  Furthermore, our result demonstrates that data about the $K_0$
structure is sufficient to distinguish up to isomorphism the algebras in an important class of purely infinite simple
$K$-algebras.
\end{abstract}

\maketitle

\section*{Introduction}

Let $K$ be any field, and let $m<n$ be positive integers. The ring
$R$ is said to have {\it invariant basis number} (IBN) if no two
free left $R$-modules of differing rank over $R$ are isomorphic.
On the other hand, $R$ is said to have {\it module type} $(m,n-m)$
in case for every pair of positive integers $a$ and $b$, (1) if
$1\leq a<m$ then the free left $R$-modules $R^a$  and $R^i$  are
not isomorphic for all positive integers $i\neq a$, and (2) if
$a,b \geq m$, then the free left $R$-modules $R^a$  and $R^b$ are
isomorphic precisely when $a\equiv b$ (mod $n-m$).  It is not hard
to show that any non-IBN ring has module type $(m,n-m)$ for some
pair of positive integers $m<n$.  (The notation used here is not
completely universal: some authors refer to the module type of
such an algebra as the pair $(m,n)$. Our notation is consistent
with that used in many of the algebra articles on this topic, and
is also consistent with the  C$^*$-algebra usage as well.)     As
shown by Leavitt in \cite{L1}, for every such pair $m,n$ there
exists a $K$-algebra $L_K(m,n)$ whose module type is $(m,n-m)$.
In particular, the module type of $L_K(1,n)$ is $(1,n-1)$. We
denote $L_K(1,n)$ by $L_n$.  Various aspects of these algebras
have been investigated, with an initial flurry of activity in the
1960's and early 1970's (e.g. \cite{B}, \cite{Co}, and \cite{L2}),
and then again in a revival beginning at the start of the new
millennium (e.g. \cite{A}, \cite{AAn1}, and \cite{AGP}).

On the ``analytic" side of the coin, Cuntz \cite{Cu} in the 1970s investigated the C$^*$-algebras $\{\mathcal{O}_n \mid
2\leq n\in \mathbb{N}\}$. There is an intimate connection between the Leavitt algebra $L_K(1,n)$ and the Cuntz algebra
$\mathcal{O}_{n}$. Specifically, for any field $K$, the elements of $L_K(1,n)$ can be viewed as linear transformations
on an infinite dimensional $K$-vector space in a natural way as a collection of shift operators. In particular, when
$K$ is the field of complex numbers,  then $L_K(1,n)$ can be viewed as acting on Hilbert space $\ell ^2$, and thereby
inherits the operator norm. The Cuntz algebra $\mathcal{O}_{n}$ is the completion of $L_{\mathbb{C}}(1,n)$ in the
metric induced by this norm.

Since $L_n \cong L_n^n$ as free left $L_n$-modules, by taking endomorphism rings we get immediately that there is a
ring isomorphism between  $L_n$ and ${\rm M}_n(L_n)$.   The first two authors extended this type of isomorphism to
additional matrix sizes in \cite{AAn1}, where they observe that $L_n \cong {\rm M}_d(L_n)$ whenever $d$ divides
$n^{\alpha}$ for some positive integer $\alpha$. In \cite{L1} Leavitt shows that for ${\rm gcd}(d,n-1)>1$, the
$K$-algebras $L_n$ and ${\rm M}_d(L_n)$ cannot be isomorphic. Since $d|n^{\alpha}$ implies ${\rm gcd}(d,n-1)=1$, these
two results yield the following natural question, posed in \cite{AAn1}, page 362:

\begin{center}
For ${\rm gcd}(d,n-1)=1$, are $L_{n}$ and ${\rm M}_d(L_{n})$ isomorphic?
\end{center}
In our main result, Theorem \ref{Main}, we answer this question in
the affirmative for all fields $K$.

Theorem \ref{Main} has important consequences in the context of C$^*$-algebras.  First, we show in Section
\ref{applications} that this result can be used to {\it directly} answer in the affirmative the following question,
posed in \cite{PS}, page 8:

\begin{center}
Are $ {\rm M}_m(\mathcal{O}_{n})$ and $\mathcal{O}_{n}$ isomorphic
whenever $m$ and $n-1$ are relatively prime?
\end{center}

While an affirmative answer to this question was provided for even $n$ in \cite{Ro}, Corollary 7.3, and subsequently
shown for all $n\geq 2$ as a consequence of \cite{Ph1}, Theorem 4.3(1),  the method we provide here is significantly
more elementary. Indeed, the second important consequence of our result is that, unlike the current situation in the
C$^*$-algebra case, the isomorphisms we present between the indicated $K$-algebras are in fact explicitly given.
Moreover, when $K=\mathbb{C}$, this explicit description carries over to an explicit description of the isomorphisms
between the appropriately sized matrix rings over Cuntz algebras.

Finally, our result demonstrates that data about the $K_0$ structure is sufficient to distinguish up to isomorphism the
algebras in an important class of purely infinite simple $K$-algebras, thus paving a path for subsequent work by the
authors \cite{AAP} towards an algebraic version of \cite{P}, Theorem 4.2.4.

The authors thank the referee for an extremely careful review of this article.

\section{Notation and basic concepts}\label{basics}

We begin by explicitly defining the Leavitt algebras $L_{K}(1,n)$.
For any positive integer $n\geq 2$, and field $K$, we denote
$L_K(1,n)$ by $L_{K,n}$, and call it the {\it Leavitt algebra of
type (1,n-1) with coefficients in $K$. } (When $K$ is understood,
we denote this algebra simply by $L_n$.)  Precisely, $L_{K,n}$ is
the quotient of the free associative $K$-algebra in $2n$
variables:
$$L_{K,n}=K<X_1,...,X_n,Y_1,...,Y_n>/T,$$
where $T$ is the ideal generated by the relations $ X_iY_j - \delta_{ij}1_K$ (for $1\leq i,j \leq n$) and
$\sum_{j=1}^{n} Y_jX_j - 1_K$. The images of $X_i,Y_i$ in $L_{K,n}$ are denoted respectively by $x_i,y_i$. In
particular, we have the equalities $x_iy_j = \delta_{ij}1_K$ and $\sum_{j=1}^{n} y_jx_j = 1_K$ in $L_n$. The algebra
$L_n$ was investigated originally by Leavitt in his seminal paper \cite{L1}.  We now list various fundamental
properties of $L_n$, culminating in the property which will serve as the focus of our investigation.

\begin{proposition}\label{Basicprops}  Let $K$ be any field.
\begin{enumerate}
\item   \cite{L1}, Theorem 8:  $L_n$ has module type $(1,n-1)$.  In particular, if $a\equiv b$ (mod $n-1$) then $L_n^a \cong
L_n^b$ as free left $L_n$-modules.  Consequently, if $a\equiv b$ (mod $n-1$), then there is an isomorphism of matrix
rings ${\rm M}_a(L_n)\cong {\rm M}_b(L_n)$.
 \item Suppose $R$ is a $K$-algebra which contains a subset
$\{a_1,...,a_n,b_1,...,b_n\}$ for which $ a_ib_j = \delta_{ij}1_R$ (for $1\leq i,j \leq n)$, and $\sum_{j=1}^{n} b_ja_j
= 1_R$. (For instance, any $K$-algebra having module type $(1,n-1)$ has this property.)   Then there exists a (unital)
$K$-algebra homomorphism from $L_n$ to $R$ extending the map $x_i\mapsto a_i$ and $y_i\mapsto b_i$ (for $1\leq i \leq
n$).
\item  \cite{L2}, Theorem 2:  $L_{n}$ is a simple $K$-algebra.
\end{enumerate}
\end{proposition}

\begin{corollary}\label{method}
Let $I$ denote the identity matrix in ${\rm M}_d(L_n)$. To show
$L_n\cong {\rm M}_d(L_n)$ it suffices to show that there is a set
$S = \{a_1,...,a_n,b_1,...,b_n\} \subseteq {\rm M}_d(L_n)$ such
that: $ a_ib_j = \delta_{ij}I$ (for $1\leq i,j \leq n$);
$\sum_{j=1}^{n} b_ja_j = I$; and $S$ generates ${\rm M}_d(L_n)$ as
a $K$-algebra.
\end{corollary}
\begin{proof}   The existence of a  nontrivial $K$-algebra homomorphism
from $L_n$ to ${\rm M}_d(L_n)$ follows from Proposition \ref{Basicprops}(2), while the injectivity of such a
homomorphism follows from Proposition \ref{Basicprops}(3). Since $\{x_1,...,x_n,y_1,...,y_n\}$ generates $L_n$ as a
$K$-algebra, the image of this homomorphism is generated by $\{a_1,...,a_n,b_1,...,b_n\} \subseteq {\rm M}_d(L_n)$.
\end{proof}

For any unital ring $R$ and $i\in\{1,2,...,d\}$ we denote the
idempotent $e_{i,i}$  of the matrix ring ${\rm M}_d(R)$ simply by
$e_i$, and we define
$$E_i = \sum_{j=1}^i e_j.$$
In this notation $E_d=I$, the identity matrix in ${\rm M}_d(R)$.

\begin{definition}\label{involutiondef}   {\rm  For any field $K$, the extension of the assignments $x_i \mapsto y_i=x_i^*$ and $y_i \mapsto x_i=y_i^*$
 for $1\leq i \leq n$ yields an involution $*$ on $L_K(1,n)$. This involution on $L_K(1,n)$ produces an involution
 on any sized matrix ring ${\rm M}_m(L_K(1,n))$ over $L_K(1,n)$ by setting $X^* =
(x_{j,i}^*)$  for each $X=(x_{i,j})\in {\rm M}_m(L_K(1,n))$.}
\end{definition}

We note that if $K$ is a field with involution (which we also denote by $*$), then a second involution on $L_K(1,n)$
may be defined
 by extending the assignments $k\mapsto k^*$ for all $k\in K$,  $x_i \mapsto y_i=x_i^*$ and $y_i \mapsto
x_i=y_i^*$ for $1\leq i\leq n$.   Of course in the case $K=\C$ we have such an involution on $K$.  Although it might be
of interest to consider this second type of involution on $L_{\C}(1,n)$ in order to maintain some natural connection
with the standard involution on the corresponding Cuntz algebra $\mathcal{O}_{n}$, we prefer to work with the
involution on $L_K(1,n)$ described in Definition \ref{involutiondef}  because it can be defined for any field $K$.  All
of the results presented in this article for involutions on $L_K(1,n)$ and their matrix rings are valid using either
type of involution.

\medskip

{\bf We now set some notation which will be used throughout the remainder of the article.}  For positive integers $d$
and $n$ we write
$$n = qd + r \mbox{ where } 1\leq r \leq d.$$
We assume throughout  that ${\rm gcd}(d,n-1)=1$, and that $d < n$. (We will relax the hypothesis $d<n$ in our main
result.) Without loss of generality we will also assume that $r\geq 2$, since $r=1$ would yield $n-1 = qd$, which along
with the hypothesis that ${\rm gcd}(d,n-1)=1$ would yield $d=1$, and the main result in this case is then the trivial
statement $L_n\cong {\rm M}_1(L_n)$. An important role will be played by the number $s$, defined as
$$s = d - (r-1).$$
Since ${\rm gcd}(d,n-1)=1$ we get also that ${\rm gcd}(s,d)=1$.

\begin{definition}\label{hsequence}
{\rm We consider the sequence $\{h_i\}_{i=1}^{d}$ of integers, whose $i^{th}$ entry is given by
$$h_i = 1 + (i-1)s \ ({\rm mod} \ d).$$
The integers $h_i$ are understood to be taken from the set $\{1,2,...,d\}$. Rephrased, we define the sequence
$\{h_i\}_{i=1}^{d}$ by setting $h_1=1,$ and, for $1\leq i \leq d-1$,
$$h_{i+1} = h_i+s \mbox{ if }h_i\leq r-1, \ \ \mbox{ and } \ \ h_{i+1}=h_i-(r-1)
\mbox{ if } h_i \geq r.$$}
\end{definition}

\smallskip

 Because ${\rm gcd}(d,s)=1$ (so that $s$ is invertible ${\rm mod} \ d$), basic number theory yields the following
\begin{lemma}\label{sequence}
$\mbox{ }$
\begin{enumerate}
\item  The entries in the sequence $h_1,h_2,...,h_d$ are distinct.
\item  The set of entries  $\{h_1,h_2,...,h_d\}$ equals the set $\{1,2,...,d\}$ (in some order).
\item  The final entry in the sequence is $r$; that is, $h_d=r$.
\end{enumerate}
\end{lemma}
\begin{proof} The only non-standard statement is (3).  Suppose $r =
1+(i-1)s \ ({\rm mod} \ d)$.  Then $r-1 = (i-1)s \ ({\rm mod} \ d)$, so that $(r-1)+s = is  \ ({\rm mod} \ d)$.  But $d
= (r-1)+s$ by definition, so this gives $d = is  \ ({\rm mod} \ d)$. Now ${\rm gcd}(s,d)=1$ gives that $i = d$, so that
$i-1=d-1$ and we get $r = 1 + (d-1)s  \ ({\rm mod} \ d) = h_d$ as desired.
\end{proof}

Our interest will lie in a decomposition of $\{1,2,...,d\}$ effected by the sequence $h_1,h_2,...,h_d$, as follows.

\begin{definition}\label{d1e1f1} {\rm We let $d_1$ denote the integer for which
$$h_{d_1} = r-1$$
in the previously defined sequence. Such an integer $d_1$ exists by  Lemma \ref{sequence}(2).  Note then that
$h_{d_1+1}= (r-1)+s = d.$   We denote by $\hat{S_1}$ the following subset of $\{1,2,...,d\}$:
$$\hat{S_1}= \{h_i | 1\leq i \leq d_1\}.$$
We denote by $\hat{S_2}$ the complement of $\hat{S_1}$ in $\{1,2,...,d\}$; in other words, $\hat{S_2} = \{h_i |
d_1+1\leq i \leq d\}$.  If we define $d_2=d-d_1$, then
$$d_1 = |\hat{S_1}|, d_2 = |\hat{S_2}|, \mbox{ and } d_1+d_2 = d.$$
Let $e_1 = |\hat{S_1} \cap \{r-1, r, r+1, ... d\}|$.  So $e_1$ is the number of elements in $\hat{S_1}$ which are at
least $r-1$. Similarly, let $e_2 = |\hat{S_2} \cap \{r-1, r, r+1, ... d\}|$. (Note by definition of $\hat{S_1}$ and
Lemma \ref{sequence}(3) we have $1,r-1 \in \hat{S_1}$ and $r,d\in \hat{S_2}$.) So we get
$$e_1+e_2=|\{r-1, r, ..., d\}| = d-(r-1)+1 = d-r+2.$$
Let $f_1 = |\hat{S_1} \cap \{1,2,...,r-1, r\}|$.  So $f_1$ is the number of elements in $\hat{S_1}$ which are at most
$r$. Similarly, let $f_2 = |\hat{S_2} \cap \{1,2,...,r\}|$. We get
$$f_1 + f_2 = r.$$   Finally, by definition we have}
$$e_1+f_1 = d_1+1 \mbox{ and } e_2 + f_2 = d_2+1.$$
\end{definition}

\begin{proposition}\label{bandt}
Write $h_{d_1}=r-1= 1 + (d_1-1)s  \ ({\rm mod} \ d)$.  So there exists a nonnegative integer $t$ with $r-1= 1 +
(d_1-1)s -td$, so that
$$r-1= 1 + (d_1-1)s -t[s+(r-1)]=1 + (d_1-1-t)s -t(r-1).$$
Let $b$ denote $d_1-1-t$.  So we have
$$r-1 = 1 + bs - t(r-1).$$
(In particular, we also have $(1+t)(r-1) = 1+bs$.) Then $e_1 = t+1$, $d_1 = 1+b+t$, and $f_1 = 1+b$.
\end{proposition}
\begin{proof} By definition, each element of the sequence $\{ h_i\}^d_{i=1}$ is the
remainder of $1+(i-1)s$ modulo $d$. Now, we will show by induction on $i$ that $h_i=1+(i-1)s-l_id$ where $l_i$ is the
number of $h_j$ for which $h_j\geq r$ and $j<i$.

For $i=1$, $h_1=1=1+(1-1)s-0d$, as $r\geq 2$ implies $l_1=0$. Now, suppose that the result holds for $i\geq 1$. If
$h_i\geq r$, $l_{i+1}=l_i+1$ by definition. Also, the computation gives us
$$h_{i+1}=h_i-(r-1)=h_i+s-d=1+(i-1)s-l_id+s-d=1+is-(l_i+1)d=1+is-l_{i+1}d.$$
On the other hand, if $h_i\leq r-1$, then $l_{i+1}=l_i$ by definition. Also, the computation gives us
$$h_{i+1}=h_i+s=1+(i-1)s-l_id+s=1+is-l_id=1+is-l_{i+1}d.$$
Thus, induction step works.

Now, for $i=d_1$, denote $l_{d_1}$ by $t$. The previous assertion shows that
$$r-1=1+(d_1-1)s-td,$$
where $t$ is the number of $h_j$ for which $h_j\geq r$ and $j<d_1$. Since $\hat{S_1}=\{ h_i\mid 1\leq i\leq d_1\}$, we
have
$$t=\vert (\hat{S_1}\setminus \{ r-1\})\cap \{ r-1, r, \dots ,d\}\vert,$$ so that $e_1=1+t$.

By definition of $b$, $d_1=1+b+t$. But $f_1=d_1+1-e_1$, so we are done.
\end{proof}

\begin{example} {\rm It will be helpful to give a specific example in order to solidify these ideas.
Suppose $n=35,d=13$.  Then ${\rm gcd}(13,35-1)=1$, so we are in the desired situation.  Now $35 = 2\cdot 13 + 9$, so
that $r=9, r-1=8,$ and $s=d-(r-1)=13-8=5$. Then the sequence $h_1,h_2,...,h_d$ is given by
$$1,6,11,3,8,13,5,10,2,7,12,4,9.$$
Since $r-1=8$, The partition $\{1,2,...,d\} = \hat{S_1}\cup \hat{S_2}$ is then
$$\{1,2,...,13\} = \{1,3,6,8,11\}\cup \{2,4,5,7,9,10,12,13\}.$$
Furthermore,
$$d_1=|\{1,3,6,8,11\}|=5, \ \ d_2=|\{2,4,5,7,9,10,12,13\}|= 8,$$
$$e_1=|\{8,11\}|=2,\ \ e_2=|\{9,10,12,13\}|=4,$$
$$f_1=|\{1,3,6,8\}|=4,\ \  f_2=|\{2,4,5,7,9\}|=5.$$
Note that $f_1 = 4 = 1+3 = 1+b$, and $e_1 = 2 = 1+1 = 1+t$. Finally, we have}
$$r-1 = 8 = 1 + 3\cdot 5 - 1\cdot 8 = 1 + bs - t(r-1).$$

\end{example}

\section{The search for appropriate matrices inside ${\rm M}_d(L_n)$}\label{findmatrices}

We start this section by giving a plausibility argument for Theorem \ref{Main}. In \cite{L1}, Theorem 5, Leavitt proves
\begin{proposition}\label{moduletypeofmatrices}
If $R$ has module type $(1,n-1)$, then ${\rm M}_d(R)$ has module type $(1, \frac{n-1}{{\rm gcd}(d,n-1)})$.
\end{proposition}
Since module type is an isomorphism invariant, this result immediately gives that $L_n$ and ${\rm M}_d(L_n)$ are not
isomorphic when ${\rm gcd}(d,n-1)
>1$.

On the other hand, in case ${\rm gcd}(d,n-1)=1$, Leavitt's proof of Proposition \ref{moduletypeofmatrices} gives an
algorithm for finding specific elements $\{a_1,...,a_n,b_1,...,b_n\}$ inside ${\rm M}_d(L_n)$ which satisfy the
appropriate relations.  So, by Corollary \ref{method}, we would be done if we could show that this set of elements
generates ${\rm M}_d(L_n)$ as a $K$-algebra.

However, this set of elements does NOT generate ${\rm M}_d(L_n)$ in general. It is instructive here to look at a
specific example. Because by \cite{AAn1}, Proposition 2.1,  we know our main result is true when $d$ divides some power
of $n$, the smallest case of interest is the situation $d=3, n=5$, since then ${\rm gcd}(d,n-1)=1$ but $d$ does not
divide any power of $n$. Leavitt's proof (for general $d,n$) manifests in this specific case that $M_3(L_5)$ has module
type $(1,4)$, and is based on an analysis of the $n=5$ elements in $M_3(L_5)$
$$X_1=\begin{pmatrix}{x_{1}}&0&0\\
    x_2&0&0\\
    x_3&0&0 \end{pmatrix}\hspace{.25in}
 X_2=\begin{pmatrix}{x_4}&0&0\\
    x_5&0&0\\
    0&x_1&0\end{pmatrix} \hspace{.25in}
 X_3=\begin{pmatrix}0&x_2&0\\
    0&x_3&0\\
    0&x_4&0\end{pmatrix} \hspace{.25in}$$

$$ X_4=\begin{pmatrix}0&x_5&0\\
    0&0&x_1\\
    0&0&x_2\end{pmatrix} \hspace{.25in}
X_5=\begin{pmatrix}0&0&x_3\\
    0&0&x_4\\
    0&0&x_5\end{pmatrix} \hspace{.25in}$$
together with the five dual matrices $Y_i = X_i^*$ for $1\leq i \leq 5$. While these ten matrices generate ``much of"
$M_3(L_5)$, these matrices do not, for instance, generate the matrix unit $e_{1,3}$.  In fact, we show below in
Proposition \ref{gradediso} that whenever ${\rm gcd}(d,n-1)=1$ but $d$ does not divide $n^{\alpha}$ for any positive
integer $\alpha$, then the matrices in ${\rm M}_d(L_n)$ which arise in the proof of \cite{L1}, Theorem 5, cannot
generate ${\rm M}_d(L_n)$.

A breakthrough in this investigation was achieved when the authors were able to show that isomorphisms between more
general structures (so-called ``Leavitt path algebras"; see e.g. \cite{AAr1}), when interpreted in light of
\cite{AAr2}, Proposition 13, in fact yield an isomorphism between $L_5$ and ${\rm M}_3(L_5)$.  By tracing through the
appropriate translation maps, the following subset of ${\rm M}_3(L_5)$ emerges as the desired set of elements, elements
which satisfy the appropriate relations {\it and} generate ${\rm M}_3(L_5)$ as a $K$-algebra:
$$X_1=\begin{pmatrix}x_{1}&0&0\\
    x_5&0&0\\
    x_3&0&0\end{pmatrix} \hspace{.25in}
 X_2=\begin{pmatrix}x_2&0&0\\
    x_4&0&0\\
    0&1&0\end{pmatrix} \hspace{.25in}
 X_3=\begin{pmatrix}0&0&{x_1}^2\\
    0&0&x_5x_1\\
    0&0&x_3x_1\end{pmatrix} \hspace{.25in}$$

$$ X_4= \begin{pmatrix}0&0&x_2x_1\\
    0&0&x_4x_1\\
    0&0&x_5\end{pmatrix} \hspace{.25in}
X_5= \begin{pmatrix}0&0&x_2\\
    0&0&x_4\\
    0&0&x_3\end{pmatrix} \hspace{.25in} $$
and $Y_i = X_i^*$ for each $1\leq i \leq 5$.  What we glean from
this particular set of matrices in ${\rm M}_3(L_5)$ is that:

\smallskip

(i) it might be useful to use $1_K$ as an entry (any number of times) in the generating matrices,

(ii)  various nonlinear monomials might play a useful role in the
generating matrices, and

(iii)  it might be of use to place elements in the matrices in
some order other than lexicographic order.

\medskip

With guidance provided by the above system of generators in ${\rm M}_3(L_5)$, one can easily check that the following
set of matrices (together with the appropriate dual matrices) is also a set of generators of ${\rm M}_3(L_5)$ which
satisfies the conditions of Corollary \ref{method}, and hence provides an isomorphism between $L_5$ and ${\rm
M}_3(L_5)$.
\begin{center}
$$X_1=\begin{pmatrix}x_{1}&0&0\\
    x_2&0&0\\
    x_3&0&0\end{pmatrix} \hspace{.25in}
 X_2=\begin{pmatrix}x_4&0&0\\
    x_5&0&0\\
    0&1&0\end{pmatrix} \hspace{.25in}$$
 $$X_3=\begin{pmatrix}0&0&{x_1}^2\\
    0&0&x_2x_1\\
    0&0&x_3x_1\end{pmatrix} \hspace{.25in}
 X_4= \begin{pmatrix}0&0&x_4x_1\\
    0&0&x_5x_1\\
    0&0&x_2\end{pmatrix} \hspace{.25in}
X_5= \begin{pmatrix}0&0&x_4\\
    0&0&x_5\\
    0&0&x_3\end{pmatrix} \hspace{.25in} $$

\end{center}

It is easy to show, and not at all unexpected,  that for each $n$, the symmetric group $S_n$ acts as automorphisms on
$L_n$ in the obvious way.  Specifically, for $\sigma \in S_n$ we define $\alpha_{\sigma}:L_n \rightarrow L_n$ by
setting $\alpha_{\sigma}(x_i)=x_{\sigma(i)}$ for each $1\leq i \leq n$, and extending linearly.   In fact, with $\sigma
\in S_5$ given by $\sigma(2)=5$, $\sigma(4)=2$, and $\sigma(5)=4$, it is straightforward to show that the corresponding
$\alpha_{\sigma}$ transforms this last set of five matrices to the previously given set.

We close this section by giving three additional sets of generating matrices for ${\rm M}_3(L_5)$.  First, consider the
set $\{X_1,X_2,X_3,X_4,X_5\}$ of matrices presented directly above.  It is relatively easy to show that by defining
$X_5'$ to be the matrix gotten by interchanging the entries $x_5$ and $x_3$ of  $X_5$, then the set
$\{X_1,X_2,X_3,X_4,X_5'\}$ (and their duals) provide a generating set for ${\rm M}_3(L_5)$.  (We note for future
reference that, in contrast, switching the entries $x_5$ and $x_4$ of $X_5$ would not provide a generating set.)

Second, consider again the set $\{X_1,X_2,X_3,X_4,X_5\}$ of matrices presented directly above. It is not difficult to
show that by defining $X_4''$ and $X_5''$ to be the matrices gotten by interchanging the entry $x_2$ of $X_4$ with the
entry $x_3$ of $X_5$, then the set $\{X_1,X_2,X_3,X_4'',X_5''\}$ (and their duals) provide a generating set for ${\rm
M}_3(L_5)$.

In Section \ref{maintheorem} we will generalize these first two observations, and show how each yields an action of
various symmetric groups as automorphisms of ${\rm M}_d(L_n)$, and hence of $L_n$, whenever ${\rm gcd}(d,n-1)=1$.

Third, and finally, it is somewhat less obvious that there are many other types of actions of various symmetric groups
on ${\rm M}_3(L_5)$.  To give one such example, here is yet another set of five matrices which, along with their duals,
provides a set of generators for ${\rm M}_3(L_5)$.  Loosely speaking, these are produced from the previous set
$\{X_1,X_2,X_3,X_4,X_5\}$ by an appropriate permutation in $S_5$ together with an interchanging of the roles of the
initial and final columns of ${\rm M}_3(L_5)$.
\begin{center}
$$X_1=\begin{pmatrix}x_{1}x_5&0&0\\
    x_2x_5&0&0\\
    x_3x_5&0&0\end{pmatrix} \hspace{.25in}
 X_2=\begin{pmatrix}x_4&0&0\\
    x_4x_5&0&0\\
    x_5^2&0&0\end{pmatrix} \hspace{.25in}$$

 $$X_3=\begin{pmatrix}x_1&0&0\\
    x_2&0&0\\
    x_3&0&0\end{pmatrix} \hspace{.25in}
 X_4= \begin{pmatrix}0&1&0\\
    0&0&x_2\\
    0&0&x_3\end{pmatrix} \hspace{.25in}
X_5= \begin{pmatrix}0&0&x_1\\
    0&0&x_4\\
    0&0&x_5\end{pmatrix} \hspace{.25in} $$

\end{center}

We will describe subsequent to the proof of Theorem \ref{Main} a number of additional, significantly different
collections of generating matrices in ${\rm M}_d(L_n)$ for ${\rm gcd}(d,n-1)=1$. Each of these collections gives rise
to an automorphism of ${\rm M}_d(L_n)$. Because Theorem \ref{Main} will demonstrate that ${\rm M}_d(L_n)\cong L_n$ for
${\rm gcd}(d,n-1)=1$, each of these automorphisms of ${\rm M}_d(L_n)$ will in turn induce an automorphism of $L_n$.

\medskip

\section{The generators of ${\rm M}_d(L_n)$}\label{generators}

In this section we present the appropriate $2n$ matrices of ${\rm M}_d(L_n)$ which generate ${\rm M}_d(L_n)$. We write
$n=qd+r$ with $2\leq r \leq d$.  We assume $d<n$, so that $q\geq 1$. The matrices $X_1, X_2, ..., X_q$ are given as
follows. For $1\leq i \leq q$ we define

$$X_i=
\begin{pmatrix}x_{(i-1)d+1}&0& &0 \\
        x_{(i-1)d+2}&0& &0 \\
        \vdots&0&...&0 \\
        x_{id}&0& &0
        \end{pmatrix}
= \sum_{j=1}^{d}x_{(i-1)d+j}e_{j,1} $$ The two matrices $X_{q+1}$
and $X_{q+2}$ play a pivotal role here.  They are defined as
follows.
$$X_{q+1}=
\begin{pmatrix}x_{qd+1}&0&0& &0&0& &0& \\
        x_{qd+2}&0&0& &0&0& &0& \\
        \vdots&0&0& &0&0& &0& \\
        x_n&0&0&...&0&0&...&0&\\
        0&1&0& &0&0& &0& \\
        0&0&1& &0&0& &0& \\
        & & &\vdots& & & & & \\
        0&0&0& ...&1&0& &0& \end{pmatrix}$$
$$ = \sum_{i=1}^{d-r}e_{i+r,i+1} + \sum_{t=1}^{r}x_{qd+t}e_{t,1}$$
and
$$X_{q+2}=
\begin{pmatrix}0& &0&1&0&0& &0&0 \\
        0& &0&0&1&0& &0&0 \\
         & & & & &\vdots& & &\\
        0& &0&0&0&0& &1&0 \\
        0&...&0&0&0&0& &0&a_{q+2,r-1} \\
        0& &0&0&0&0& &0&a_{q+2,r} \\
         & & &\vdots& & & & &\vdots \\
        0& &0&0&0&0& &0&a_{q+2,d}\end{pmatrix}$$
$$ = \sum_{j=1}^{r-2}e_{j,j+s} + \sum_{t=1}^{d-(r-2)}a_{q+2,(r-2)+t}e_{(r-2)+t,d}$$
(where the elements $a_{q+2,r-1}, a_{q+2,r},..., a_{q+2,d} \in
L_n$ are monomials in $x$-variables which will be determined
later).  In case $d-r=0$ or $r-2=0$ we interpret the appropriate
sums as zero.

The remaining matrices $X_{q+3},...,X_n$ will be explicitly
specified later, but each of these will have the same general
form. In particular, for $q+3 \leq i \leq n$,
$$X_i=
\begin{pmatrix}0& &0&a_{i,1} \\
        0& &0&a_{i,2} \\
         0&...&\vdots&\\
        0& &0&a_{i,d}
        \end{pmatrix}
= \sum_{j=1}^{d}a_{i,j}e_{j,d} $$ (where the elements $a_{i,1},
a_{i,2}, ... ,a_{i,d} \in L_n$ are monomials in the $x$-variables
which will be determined later).  In case $q+3 > n$ then we
understand that there are no matrices of this latter form in our
set of $2n$ matrices.  We note that we always have the matrices
$X_{q+1}$ and $X_{q+2}$, since $n=qd+r\geq q\cdot 1 + 2$.

\bigskip

We define the matrices $Y_i$ for $1\leq i \leq n$ by setting $Y_i
= X_i^*$. Because they will play such an important role, we
explicitly describe $Y_{q+1}$ and $Y_{q+2}$.

$$Y_{q+1}=
\begin{pmatrix}y_{qd+1}&y_{qd+2}&...&y_n&0&0& &0& \\
        0&0&0&0&1&0& &0& \\
        0&0&0&0&0&1& &0& \\
         & & &\vdots& & &...& &\\
        0&0&0&0&0&0& &1& \\
        0&0&0&0&0&0& &0& \\
        & & &\vdots& & & & & \\
        0&0&0&0&0&0& &0& \end{pmatrix}$$
$$ = \sum_{i=1}^{d-r}e_{i+1,i+r} + \sum_{t=1}^{r}y_{qd+t}e_{1,t}$$
and
$$Y_{q+2}=
\begin{pmatrix}0&0& &0&0&0& &0 \\
         & & & &\vdots& & &  \\
        0&0& &0&0&0& &0\\
        1&0& &0&0&0&... &0 \\
        0&1& &0&0&0& &0 \cr
         & &...& & & & &\vdots \\
        0&0& &1&0&0& &0 \\
        0&0& &0&a^*_{q+2,r-1}& a^*_{q+2,r}&...& a^*_{q+2,d} \end{pmatrix}$$
$$ = \sum_{j=1}^{r-2}e_{j+s,j} + \sum_{t=1}^{d-(r-2)}a^*_{q+2,(r-2)+t}e_{d,(r-2)+t}$$
As above, in case $d-r=0$ or $r-2=0$ we interpret the corresponding sums as zero.

\begin{definition}
{\rm We denote by $A$ the subalgebra of ${\rm M}_d(L_n)$ generated by the matrices $$\{X_i,Y_i | 1\leq i \leq n\}.$$
That is,}
$$A =<\{X_i,Y_i | 1\leq i \leq n\}>.$$
\end{definition}

So in order to achieve our main result, we seek to show that $A={\rm M}_d(L_n)$.

\medskip

Using the relation $\sum_{j=1}^n y_jx_j = 1_K$, we immediately get

\begin{lemma}\label{YiXi1toqplus1}
$$\sum_{i=1}^{q+1}Y_iX_i = E_s \in A.$$
\end{lemma}

A similar computation yields

\begin{lemma}\label{IminusEsinA}
$\mbox{ }$
\begin{enumerate}
\item Assume the elements $\{a_{q+2,r-1},...,a_{q+2,d}\} \cup \{a_{i,j} | q+3 \leq i \leq n, 1\leq j \leq d\}$ are chosen
so that
$$
\hspace{6.5truecm} \sum_{i,j}a^*_{i,j}a_{i,j}=1_K. \hspace{6truecm}
(\dagger)
$$
Then
$$\sum_{i=q+2}^{n}Y_iX_i = I - E_s\in A.$$
\item Assume the elements $\{a_{q+2,r-1},...,a_{q+2,d}\}$ are chosen
so that
$$a_{q+2, j}a^*_{ q+2, i}=\delta _{i,j}$$
for every $i,j\in \{r-1, \dots , d\}.$ Then $X_{q+2}Y_{q+2}=I$.
\end{enumerate}
\end{lemma}

\begin{lemma}\label{dual}
If $a\in A$ then $a^*\in A$.
\end{lemma}
\begin{proof} The set of generators of $A$ has this property, and the relations
are self-dual, hence for any element $a$ which can be generated by
ring-theoretic operations we can also generate $a^*$.
\end{proof}

\begin{definition}
{\rm Recall the partition $\hat{S_1} \cup \hat{S_2}$ of $\{1,2,...,d\}$ described in Section \ref{basics}. For $i,j\in
\{1,2,...,d\}$ we write $i\sim j$ in case $i,j$ are both in the same $\hat{S_k}, k=1,2$.}
\end{definition}

Our goal for the remainder of this section is to show that $A$
contains all matrix units $e_{i,j}$ for $i\sim j$. We begin by
defining two monomorphisms of ${\rm M}_d(L_n)$ which will be useful
in this context.

\begin{definition}
{\rm We define the monomorphism $\beta$ of ${\rm M}_d(L_n)$ by setting
$$\beta(M) = Y_{q+1}MX_{q+1}$$
for each $M\in {\rm M}_d(L_n).$   Since $Y_{q+1}$ and $X_{q+1}$ are each in $A$, then $\beta$ in fact restricts to a
monomorphism of $A$.

Assuming that we have chosen the elements
$\{a_{q+2,r-1},...,a_{q+2,d}\}$ as described in Lemma
\ref{IminusEsinA}(2), we define the monomorphism $\phi$ of ${\rm
M}_d(L_n)$ by setting
$$\phi(M) = Y_{q+2}MX_{q+2}$$
for each $M\in {\rm M}_d(L_n).$   Since $Y_{q+2}$ and $X_{q+2}$ are each in $A$, then $\phi$ in fact restricts to a
monomorphism of $A$.}

\end{definition}

We begin by showing that all of the matrix idempotents
$\{e_i|1\leq i\leq d\}$ are in $A$. The results presented in the
next two lemmas follow directly from straightforward matrix
computations, so we omit their proofs.

\begin{lemma}\label{klessthanrminus1}
If $k<r-1$ then
$$\phi(e_k) = e_{k+s}.$$
\end{lemma}

\begin{lemma}\label{kbiggerthanr}
If $k>r$ then
$$\beta(e_k) = e_{k-(r-1)}.$$
\end{lemma}

It is instructive to note the following. In words, the previous two lemmas say that we can move matrix idempotents
``forward by $s$" (if we start with an index less than $r-1$), and ``backwards by $r-1$" (if we start with an index
bigger than $r$). But even though it would make sense to move the specific idempotent $e_{r-1}$ forward by $s$ units
(since $(r-1)+s=d$), or to move the specific idempotent $e_{r}$ backwards by $r-1$ units, neither of these moves can be
effected by the matrix multiplications described in the lemmas.  For instance, the entry in the $(d,d)$ coordinate of
$\phi(e_{r-1}) = Y_{q+2}e_{r-1}X_{q+2}$ is $a^*_{q+2,r-1}a_{q+2,r-1}$, which may or may not equal $1$ depending on the
choice of $a_{q+2,r-1}$.  (Indeed, we will see later that we will NOT choose $a_{q+2,r-1}$ having this property.) This
observation is precisely the reason why we must expend so much effort in analyzing the partition $\hat{S_1}\cup
\hat{S_2}$ of $\{1,2,...,d\}$ described previously.

\medskip

We consider the sequence $\{u_i\}_{i=1}^{d}$ of integers, whose
$i^{th}$ entry is given by
$$u_i = is \ ({\rm mod} \ d).$$
The integers $u_i$ are understood to be taken from the set
$\{1,2,...,d\}$.  Rephrased, we define the sequence
$\{u_i\}_{i=1}^{d}$ by setting $u_1=s,$ and, for $1\leq i \leq
d-1$,
$$u_{i+1} = u_i+s \mbox{ if }u_i\leq r-1, \ \ \mbox{ and } \ \ u_{i+1}=u_i-(r-1)
\mbox{ if } u_i \geq r.$$ Of course, the $u$-sequence is closely related to the $h$-sequence described in Section
\ref{basics}.  Thus it is not surprising that the following Lemma closely resembles Lemma \ref{sequence}.  Because
${\rm gcd}(d,s)=1$ (so that $s$ is invertible ${\rm mod} \ d$), basic number theory yields the following

\begin{lemma}\label{usequence}
$\mbox{ }$
\begin{enumerate}
\item   The entries in the sequence $u_1,u_2,...,u_d$ are distinct.
\item  The set of entries  $\{u_1,u_2,...,u_d\}$ equals the set $\{1,2,...,d\}$ (in some order).
\item  The penultimate entry in the sequence is $r-1$; that is, $u_{d-1}=r-1$.
\item  The final entry in the sequence is $d$; that is, $u_d=d$.
\end{enumerate}
\end{lemma}
\begin{proof} The only non-standard statements are (3) and (4).
Suppose $r-1 = is  \ ({\rm mod} \ d)$.  Then $d = r-1 + s = (i+1)s \ ({\rm mod} \ d)$.  Now ${\rm gcd}(s,d)=1$ gives
that $i+1 = d$, so that $i=d-1$ and we get $r-1=(d-1)s  \ ({\rm mod} \ d) = u_{d-1}$ as desired.   Then (4) follows
directly from (3) and the equation $d=(r-1)+s$.
\end{proof}

\begin{proposition}\label{eiinA}
For every $j$ with $1\leq j \leq d$ we have $e_j\in A$.
\end{proposition}
\begin{proof} The key idea is to show that $E_j\in A$ for all $1\leq j
\leq d$. Since $X_1Y_1 = I$ we have $I = E_d \in A$.  We consider the sequence of matrices
$E_{u_1},E_{u_2},...,E_{u_d}$ arising from the sequence $\{u_i\}_{i=1}^{d}$ described above. By induction on $i$, we
show that each of $E_{u_1},E_{u_2},...,E_{u_{d-1}}\in A$. For $i=1$ we have $E_{u_1}=E_s\in A$ by Lemma
\ref{YiXi1toqplus1}. Now we assume that $E_{u_i}\in A$ for $i\leq d-2$, and show that $E_{u_{i+1}}\in A$. By Lemma
\ref{usequence}(3), $i\leq d-2$ gives that $u_i\neq r-1$. There are two cases.

\smallskip

Case 1:  $u_i \leq r-2$. Then by definition $u_{i+1}=u_i+s$. Since
$E_{u_i}\in A$ by hypothesis, we have $\phi(E_{u_i})\in A$, which
then gives
$$E_s + \phi(E_{u_i})\in A.$$
But since $u_i\leq r-2$, Lemma \ref{klessthanrminus1} applies to
give
$$\phi(E_{u_i}) = \sum_{j=1}^{u_i} e_{j+s} = \sum_{k=s+1}^{u_i+s} e_{k},$$
so that
$$E_s+\phi(E_{u_i})=\sum_{k=1}^{s} e_{k} + \sum_{k=s+1}^{u_i+s} e_{k}= E_{u_i+s}=E_{u_{i+1}},$$
so that $ E_{u_{i+1}}\in A$, and Case 1 is shown.

\smallskip

Case 2:  $u_i\geq r$. Since $I = E_d \in A$ we have $I-E_{u_i} \in
A$, and since $E_s\in A$ we get
$$E_s - \beta(I-E_{u_i} )\in A.$$
But $I-E_{u_i} = \sum_{j=u_i+1}^d e_j$, and $u_i+1 > r$, so Lemma
\ref{kbiggerthanr} applies to give
$$\beta(I-E_{u_i}) = \sum_{j=u_i+1}^d e_{j-(r-1)}.$$
Thus we get that
\begin{eqnarray*}
E_s - \beta(I-E_{u_i}) & = & E_s - \sum_{j=u_i+1}^d e_{j-(r-1)} \\
& = &  E_s - \sum_{k=u_i+1-(r-1)}^{d-(r-1)} e_{k} = E_s - \sum_{k=u_i+1-(r-1)}^{s} e_{k} \\
& = &  \sum_{k=1}^{u_i-(r-1)} e_{k}=E_{u_i-(r-1)} = E_{u_{i+1}},
\end{eqnarray*}
so that $E_{u_{i+1}}\in A$, and Case 2 is shown. Thus we have established by induction that $E_{u_i}\in A$ for all
$1\leq i \leq d-1$.   But $E_{u_d}=E_d$ by Lemma \ref{usequence}, and $E_d =I\in A$ has already been established, so in
fact we have $E_{u_i}\in A$ for all $1\leq i \leq d.$  So by Lemma \ref{usequence}(2) we conclude that $E_j\in A$ for
all $1\leq j\leq d$.

Now the desired result follows easily from the observation that
$e_1=E_1\in A$, while $e_j=E_j - E_{j-1} \in A$ for all $2\leq
j\leq d$.
\end{proof}

We remark that we need not modify the proof of Proposition \ref{eiinA} at all in case $r=2$ (resp. $r=d$). This is
because even though we would not have the matrix $X_{q+2}$ (resp. $X_{q+1}$) containing 1 in the appropriate entries,
in the case $r=2$ (resp. $r=d$) we would have $s=d-1$ (resp. $s=1$), so that we would only be using multiplication by
$X_{q+1}$ (resp. $X_{q+2}$) in the proof.

\smallskip

Now that we have established that all of the matrix idempotents
$e_i$ ($1\leq i \leq d$) are in $A$, we use them to generate all
of the matrix units $e_{i,j}$.

\begin{lemma}\label{e11plussinA}
$\mbox{ }$
\begin{enumerate}
\item  Suppose $1+s < d$.  Then $e_1X_{q+2}e_{1+s} = e_{1,1+s}$, and $e_{1,1+s}\in A$.
\item  Suppose $1+s =d$.  Then $\hat{S_1}=\{1\}$ and $\hat{S_2}=\{2,...,d\}$.
\item  The situation  $1+s >d$ is not possible.
\end{enumerate}
\end{lemma}
\begin{proof}
(1) We have $1+s\leq d-1$. By construction, the $(1,s+1)$ entry of $X_{q+2}$ is $1$ as long as $r-2\geq 1$. But
$1+s\leq d-1$ gives $d-(r-2)\leq d-1$, which yields the desired $r-2\geq 1$. Now use Proposition \ref{eiinA}.

(2) If $1+s =d$, since $(r-1)+s=d$ we get $r-1=1$.  So the sequence $\{h_i\}_{i=1}^d$ has $h_1=1=r-1$, so that
$\hat{S_1}=\{1\}$.

(3) If $1+s >d$, then  with $(r-1)+s=d$ we would get  $r< 2$, contradicting the hypothesis that $r\geq 2$.
\end{proof}

\begin{lemma}\label{edsinA}
$\mbox{ }$
\begin{enumerate}
\item  Suppose $n$ is not a multiple of $d$. Then $e_dX_{q+1}e_{s} = e_{d,s}$, and $e_{d,s}\in A$.
\item  Suppose $n$ is a multiple of $d$. Then
$\hat{S_1}=\{1,2,...,d-1\}$ and $\hat{S_2}=\{d\}$.
\end{enumerate}
\end{lemma}
\begin{proof}
(1) If $n$ is not a multiple of $d$ then $r\neq d$, so that the $(d,s)$ entry of the matrix  $X_{q+1}$ is $1$.  Now use
Proposition \ref{eiinA}.

(2) On the other hand, if $n$ is a multiple of $d$, then $n=qd+d$, so $r=d$, so that $r-1=d-1$, which gives
$s=d-(r-1)=1$, so that the sequence $\{h_i\}_{i=1}^d$ has $h_1=1$, $h_i=h_i+1$, and the result follows.
\end{proof}

The next Proposition provides a link between the matrix units $e_{i,j}\in A$ and the partition $\hat{S_1}\cup
\hat{S_2}$ of $\{1,2,...,d\}$.

\begin{proposition}\label{partitionconsecutivematrixunits}
Consider the sequence $\{h_i\}_{i=1}^d$ described in Section
\ref{basics}. Let $h_i,h_{i+1},h_{i+2}$ be three consecutive
elements of the sequence, where $h_i\neq r,r-1$ and $h_{i+1}\neq
r,r-1$.  (In other words, consider three consecutive elements
$h_i,h_{i+1},h_{i+2}$ so that all three are in $\hat{S_1}$ or all
three are in $\hat{S_2}$.) Then there exists $X\in
\{X_{q+1},X_{q+2}\}$ and $Y\in \{Y_{q+1},Y_{q+2}\}$ so that
$$Ye_{h_i,h_{i+1}}X = e_{h_{i+1},h_{i+2}}.$$
In particular, in this situation, if $e_{h_i,h_{i+1}}\in A$ then
also $e_{h_{i+1},h_{i+2}}\in A$.
\end{proposition}
\begin{proof}
There are four cases to consider, depending on whether we
use the ``plus $s$" or ``minus $r-1$" operation to get from one
element of the sequence to the next.

\medskip

Case 1:  $h_{i+1}=h_i+s$ and $h_{i+2}=h_{i+1}+s$.  In this
situation we have $h_i\leq r-1$ because $h_{i+1}\leq d=s+(r-1)$
and $h_{i+1}=h_i+s$.  But $h_i \neq r-1$ by hypothesis.  Thus we
have in fact $h_i\leq r-2$.   In an exactly analogous way we also
have $h_{i+1}\leq r-2$.   Using that each of $h_i$ and $h_{i+1}$
is less than $r-1$, we get
$$Y_{q+2}e_{h_i,h_{i+1}}X_{q+2} = e_{h_{i+1},h_{i+2}}.$$

\smallskip

Case 2:  $h_{i+1}=h_i+s$ and $h_{i+2}=h_{i+1}-(r-1)$.  As in Case
1 we have $h_i< r-1$.  Also, $h_{i+1}\geq r$ because $1\leq
h_{i+2}=h_{i+1}-(r-1)$.  But $h_{i+1}\neq r$ by hypothesis.  Thus
we have in fact $h_{i+1}>r$.  Using both that $h_i<r-1$ and
$h_{i+1}>r$, we get
$$Y_{q+2}e_{h_i,h_{i+1}}X_{q+1} = e_{h_{i+1},h_{i+2}}.$$

\smallskip

Case 3:  $h_{i+1}=h_i-(r-1)$ and $h_{i+2}=h_{i+1}+s$. As shown
above, the hypotheses yield $h_i>r$ and $h_{i+1}<r-1$, from which
we get
$$Y_{q+1}e_{h_i,h_{i+1}}X_{q+2} = e_{h_{i+1},h_{i+2}}.$$

\smallskip

Case 4:  $h_{i+1}=h_i-(r-1)$ and $h_{i+2}=h_{i+1}-(r-1)$. As shown
above, the hypotheses yield $h_i>r$ and $h_{i+1}>r$, from which we
get
$$Y_{q+1}e_{h_i,h_{i+1}}X_{q+1} = e_{h_{i+1},h_{i+2}},$$
and the result is established.
\end{proof}

We now establish the relationship between the partition $\hat{S_1}\cup \hat{S_2}$ of $\{1,2,...,d\}$ and the matrix
units $e_{i,j}\in A.$  Intuitively, the idea is this.  Suppose for instance that $a,b\in \hat{S_1}$. We seek to show
that $e_{a,b}\in A$.  There is a sequence of elements in $\hat{S_1}$ which starts at $a$ (resp. $b$) and ends at $r-1$.
By the previous result, this will imply that $e_{a,r-1}\in A$ (resp. $e_{b,r-1}\in A$).  But then by duality
$e_{r-1,b}\in A$, so that $e_{a,r-1}e_{r-1,b}=e_{a,b}\in A$.   Here are the formal details.

\begin{proposition}\label{partitionmatrixunits}
Suppose $h_i,h_{j}$ are two entries in the sequence
$\{h_i\}_{i=1}^d$, for which both entries are either in
$\hat{S_1}$ or $\hat{S_2}$.  Then $$e_{h_i,h_j}\in A.$$
\end{proposition}
\begin{proof} We start by proving the result for $\hat{S_1}$. Suppose
first that we are in a situation for which $1+s < d$. Then Lemma \ref{e11plussinA}(1) yields that $e_{1,1+s}\in A$.
Since in this situation the integers $1,1+s$ are the first two elements of the sequence $\{h_i\}_{i=1}^d$, and both are
in $\hat{S_1}$, repeated applications of Proposition \ref{partitionconsecutivematrixunits} gives that
$e_{h_i,h_{i+1}}\in A$ for any two consecutive elements $h_i,h_{i+1}$ of $\hat{S_1}$.  By matrix multiplication this
then gives $e_{h_i,h_{j}}\in A$ whenever $i<j$ and both $h_i,h_{j}$ are in $\hat{S_1}$.  By Lemma \ref{dual} this gives
that $e_{h_i,h_{j}}\in A$ whenever $i\neq j$ and both $h_i,h_{j}$ are in $\hat{S_1}$.  This together with Proposition
\ref{eiinA} yields that $e_{h_i,h_{j}}\in A$ whenever both $h_i,h_{j}$ are in $\hat{S_1}$.

On the other hand, if we are in a situation for which $1+s = d$, then by Lemma \ref{e11plussinA}(2) we have that
$\hat{S_1}=\{1\}$, and the result follows immediately from Proposition \ref{eiinA}.

The result for $\hat{S_2}$ is established in a similar manner,
using Lemma \ref{edsinA} and Propositions \ref{eiinA} and
\ref{partitionconsecutivematrixunits}, along with the fact that
whenever $n$ is not a multiple of $d$, then the first two elements
of $\hat{S_2}$ in the sequence $\{h_i\}_{i=1}^d$ are $d,s$.
\end{proof}

\section{The main theorem}\label{maintheorem}

With the results of Section \ref{generators} in hand, we now show
how the partition $\hat{S_1}\cup \hat{S_2}$ of $\{1,2,...,d\}$ can
be used to specify the elements of $X_{q+2},...,X_n$ in such a way
that the set $$\{X_1,...,X_n,Y_1,...,Y_n\}$$ generates ${\rm
M}_d(L_n)$.

\begin{definition}
{\rm We define a partition $S_1 \cup S_2$ of $\{1,2,...,n\}$ as follows: For $w\in \{1,2,...,n\}$, write $w = q_wd +
\hat{w}$ with $1\leq \hat{w} \leq d$. We then define  $w \in S_k$ (for $k= 1,2$) if and only if $\hat{w}\in
\hat{S_k}$.}
\end{definition}

So we are 'enlarging' the partition of $\{1,2,...,d\} = \hat{S_1}
\cup \hat{S_2}$ to a partition of $\{1,2,...,n\} = S_1 \cup S_2$
by extending modulo $d$.


Now consider this set, which we will call ``The List":

$$x_1^{d-1}$$
$$x_2x_1^{d-2}, x_3x_1^{d-2}, ... , x_nx_1^{d-2}$$
$$x_2x_1^{d-3}, x_3x_1^{d-3}, ... , x_nx_1^{d-3}$$
$$\vdots$$
$$x_2x_1, x_3x_1, ... , x_nx_1$$
$$x_2, x_3, ... , x_n$$

\begin{lemma}\label{Listhasdagger}  The elements of The List satisfy $(\dagger)$. That is,

$$ y_1^{d-1}x_1^{d-1} + \sum_{i=0}^{d-2} \sum_{j=2}^n y_1^iy_j x_jx_1^i = 1_K.$$

\end{lemma}
\begin{proof} We note that
\begin{eqnarray*}
 y_1^{d-1}x_1^{d-1} + \sum_{i=0}^{d-2} \sum_{j=2}^n y_1^iy_j x_jx_1^i & = &
y_1^{d-1}x_1^{d-1} + \sum_{j=2}^n y_1^{d-2}y_j x_jx_1^{d-2} + \sum_{i=0}^{d-3} \sum_{j=2}^n y_1^iy_j x_jx_1^i \\
 & = &  \sum_{j=1}^n y_1^{d-2}y_j x_jx_1^{d-2} + \sum_{i=0}^{d-3} \sum_{j=2}^n y_1^iy_j x_jx_1^i \\
 & = & y_1^{d-2} (\sum_{j=1}^n y_j x_j) x_1^{d-2} + \sum_{i=0}^{d-3} \sum_{j=2}^n y_1^iy_j x_jx_1^i \\
& = &  y_1^{d-2}(1_K) x_1^{d-2} + \sum_{i=0}^{d-3} \sum_{j=2}^n y_1^iy_j x_jx_1^i \\
& = &  y_1^{d-2} x_1^{d-2} + \sum_{i=0}^{d-3} \sum_{j=2}^n y_1^iy_j x_jx_1^i.
\end{eqnarray*}

By induction we continue in a similar way to get
$$  \ \ \ =  y_1 x_1 +  \sum_{j=2}^n y_j x_j = 1_K. $$
\end{proof}

It is clear that

\begin{lemma}\label{listsize}
There are $(d-1)(n-1)+1$ elements on The List.
\end{lemma}

\begin{lemma}\label{numbertobespecified}  The number of entries
$\{a_{q+2,r-1},...,a_{q+2,d}\}\cup \{a_{i,j} | q+3 \leq i \leq n, 1\leq j \leq d\}$ which must be specified to form the
matrices $X_{q+2}, X_{q+3}, ..., X_n$ is
$$(s+1)+d[n-(q+2)].$$
\end{lemma}
\begin{proof} The elements $\{a_{q+2,j} | r-1 \leq j \leq d\}$ needed to
complete $X_{q+2}$ is a list containing $d-(r-1)+1 = s+1$ entries. There are $n-(q+2)$ matrices in the list $X_{q+3},
...,X_n$, and each of these matrices will contain exactly $d$ nonzero entries.
\end{proof}

\begin{lemma}
The number of entries which must be specified to form the matrices $X_{q+2}, X_{q+3}, ...,X_n$ is equal to the number
of entries in The List.
\end{lemma}
\begin{proof} By Lemmas \ref{listsize} and \ref{numbertobespecified}, we
must show
$$(s+1)+d[n-(q+2)] = (d-1)(n-1)+1.$$
But \begin{eqnarray*}
 (s+1)+d[n-(q+2)] & = & [d-(r-1)]+1 + dn-dq-2d \\
 & = &  d-r+2+dn-dq-2d \\
 & = &  -(n-qd) +2+dn-dq-d \\
 & = &  n(d-1) +2 -d \\
 & = &  n(d-1) -(d-1) + 1 \\
 & = &  (n-1)(d-1)+1.
 \end{eqnarray*}
\end{proof}

The following result describes exactly how many of the entries to
be specified in $X_{q+2},...,X_n$ correspond to the subset
$\hat{S_1}$ in the partition $\hat{S_1}\cup \hat{S_2}$ of
$\{1,2,...,d\}$.

\begin{lemma}\label{specifiedinS1hat}
Consider the set of matrices $X_{q+3}, ..., X_n$, together with the last $s+1$ rows of $X_{q+2}$.  Then the number of
nonzero entries corresponding to rows indexed by elements of $\hat{S_1}$ equals
$$d_1[n-(q+2)]+e_1.$$
\end{lemma}
\begin{proof} This follows directly by an argument analogous to that given
in the proof of Lemma \ref{numbertobespecified}, together with the
definitions of $d_1$ and $e_1$.
\end{proof}

\begin{lemma}\label{thelistinS1}
The number of entries on The List of the form $x_ux_1^t$ for which
$u\in S_1$ is
$$(d-1)[(qd_1-1)+f_1]+1.$$
\end{lemma}
\begin{proof} Consider each of the $d-1$ rows of The List (other than the
first).  For each of the $d_1$ entries which are in $\hat{S_1}$
(including $1$) there are $q$ elements congruent to it (modulo
$d$).  So we get $qd_1$ such entries. But we have started each
list with $x_2$ (and not $x_1$), so in fact there are $qd_1 - 1$
such entries in each row. Each row also contains $f_1$ entries
from the set $\{qd+1, ..., qd+r=n\}$.   There are $d-1$ rows.
Finally, we add in the term corresponding to $x_1^{d-1}$.
\end{proof}

Before we get to the main proposition, we need a computational
lemma.

\begin{lemma}\label{d1r}
$$d_1r = df_1 - d + d_1 + 1.$$
\end{lemma}
\begin{proof} Using the equations $d_1=1+b+t$ and $(1+t)(r-1) = 1+bs$ from
Proposition \ref{bandt}, we get
 \begin{eqnarray*}
 d_1r & = & (1+b+t)r = (1+t)r+ br \\
      & = &(1+t)(r-1) + (1+t) + br \\
      & = & bs+br+t+2
 \end{eqnarray*}
  while
\begin{eqnarray*}
df_1-d+d_1+1 & = &(s+(r-1))(1+b)-(s+(r-1))+(1+b+t)+1 \\
             & = & sb+rb+t+2 \ \ \mbox{ (by an easy computation). }
\end{eqnarray*}
\end{proof}

We are now ready to prove the key algorithmic tool which will provide the vehicle for our main result.

\begin{proposition}\label{possibleplacement}
Consider the set of matrices $X_{q+3}, ..., X_n$, together with the last $s+1$ rows of $X_{q+2}$.  Then the number of
nonzero entries corresponding to rows indexed by elements of $\hat{S_1}$ equals the number of entries on The List of
the form $x_ux_1^t$ for which $u\in S_1$.

Rephrased:  It is possible to place the elements of The List in the ``to be specified" entries of the matrices
$X_{q+2},X_{q+3}, ... ,X_n$ in such a way that each entry of the form $x_ux_1^t$ for $u\in S_k$ ($k=1,2$) is placed in
a row indexed by $\hat{u}$ where $\hat{u}\in \hat{S_k}$ ($k=1,2$).
\end{proposition}
\begin{proof} By Lemmas \ref{specifiedinS1hat} and \ref{thelistinS1} it
suffices to show that
$$d_1[n-(q+2)]+e_1 = (d-1)[(qd_1-1)+f_1]+1.$$
But
 \begin{eqnarray*}
d_1[n-(q+2)] + e_1  & = & d_1[qd+r-q-2]+e_1 \\
        & = & d_1q(d-1) + d_1r - 2d_1 + e_1 \\
        & = & d_1q(d-1) + [df_1 -d + d_1 + 1] -2d_1 + e_1 \mbox{ (using Lemma \ref{d1r}) } \\
        & = & d_1q(d-1) + [df_1 -d + d_1 + 1] -2d_1 + [d_1+1-f_1] \mbox{ (Definition \ref{d1e1f1}) } \\
        & = & d_1q(d-1) + (d-1)f_1 - (d-1) + 1 \\
        & = & (d-1)[d_1q -1 + f_1]  + 1
\end{eqnarray*}
and we are done.
\end{proof}

In other words, Proposition \ref{possibleplacement} implies that is possible to place the entries of The List in the
empty ``boxes" of the matrices $X_{q+2},X_{q+3}, ..., X_n$ in such a way that each entry of the form $x_ux_1^t$ for
$u\in S_k$ ($k=1,2$) is placed in a row indexed by $\hat{u}$ where $\hat{u}\in \hat{S_k}$ ($k=1,2$).

\medskip

{\bf We assume for the remainder of this article that we have made such a placement.} To help the reader clarify the
process, a specific example appears below.  However, the reader should keep in mind that in fact there are {\it many}
possible such placements.

\medskip

Once such a placement has been made, we can immediately deduce various properties of the matrices
$\{X_1,...,X_n,Y_1,...,Y_n\}$.   For instance,

\begin{lemma}\label{XYequalsI}
  For all $1\leq i,j \leq n$ we have
$$X_iY_j = \delta_{i,j}I \mbox{ in } {\rm M}_d(L_n).$$
\end{lemma}
\begin{proof}
By definition of the matrices $X_i,Y_j$ it suffices to show that
$$x_ix_1^t \cdot y_1^u y_j = \delta_{t,u}\delta_{i,j}1_K$$
for all $1\leq i,j\leq n$ and $1\leq u,t \leq d-2$.  But this follows easily by the definition of multiplication in
$L_n$.
\end{proof}

\begin{lemma}\label{eachw}
For each $w$ having $1\leq w \leq n$, $x_we_{\hat{w},1}\in A$ where $w\sim \hat{w}$.
\end{lemma}
\begin{proof}  Write $w = q_wd+ \hat{w}$ with $1\leq \hat{w}\leq d$. But
$e_{\hat{w}}$ and $e_1$ are in $A$, so $e_{\hat{w}}X_{q_w}e_1 \in A$, and this gives the result.
\end{proof}

\begin{lemma}\label{eachv}
For each $v$ having $1\leq v \leq n$, $y_ve_{1,\hat{v}}\in A$ where $v\sim \hat{v}$.
\end{lemma}
\begin{proof}  Write $v = q_vd+ \hat{v}$ with $1\leq \hat{v} \leq d$.  But
$e_{\hat{v}}$ and $e_1$ are in $A$, so $e_1Y_{q_v}e_{\hat{v}} \in A$, and this gives the result.
\end{proof}
\noindent
(We note that indeed Lemma \ref{eachv} can also be established directly from Lemmas \ref{eachw} and
\ref{dual}.)

\medskip

Proposition  \ref{partitionmatrixunits} yields that matrix units indexed by the sets $\hat{S_1}$ and $\hat{S_2}$ are in
$A$. In order to show that all the matrix units $\{e_{i,j}|1\leq i,j\leq d\}$ are in $A$, we need to provide a
``bridge" between these two subsets of matrix units. That connection is made in the following Proposition, which
provides the last major piece of the puzzle.

\begin{proposition}\label{e1dinA}
$$e_{1,d}\in A \ \mbox{ and } \ e_{d,1}\in A.$$
\end{proposition}
\begin{proof}
Because we have assumed that we have placed the elements from The List in a manner ensured by Proposition
\ref{possibleplacement}, there exists $M \in \{X_{q+2},X_{q+3},...,X_n\}$ and an integer $l \in \{1,2,...,d\}$ for
which $l\sim 1$, and for which the $(l,d)$ entry of $M$ is $x_1^{d-1}$.  That is,  $e_{l}Me_d= x_1^{d-1}e_{l,d}\in A$.
But because $l\sim 1$, Proposition \ref{partitionmatrixunits} gives that $e_{1,l}\in A$.  Thus
$e_{1,l}x_1^{d-1}e_{l,d}\in A$, so
$$x_1^{d-1}e_{1,d}\in A.$$
We have $y_1e_1 = e_1Y_1e_1\in A$, so that
$$y_1e_1\cdot x_1^{d-1}e_{1,d}=y_1x_1^{d-1}e_{1,d} =
y_1x_1x_1^{d-2}e_{1,d}\in A.$$

Now choose any $w$ with $2\leq w \leq n$.  Again using the
hypothesis that we have placed the elements from The List in a
manner ensured by Proposition \ref{possibleplacement}, there exists
$M \in \{X_{q+2},X_{q+3},...,X_n\}$ and $w'\in \{1,...,d\}$ for
which $w'\sim w$, and
$$e_{w'}Me_d=x_wx_1^{d-2}e_{w',d}\in A.$$

 Write $w = q_wd + \hat{w}$ with $1\leq \hat{w}\leq d$. Then
$w\sim \hat{w}$ by definition, and so we get $w'\sim \hat{w}$. So by
Proposition \ref{partitionmatrixunits}, $e_{\hat{w},w'}\in A$. In
addition, $e_1Y_{q_w}e_{\hat{w}}=y_we_{1,\hat{w}}\in A$.     So we
get
$$y_we_{1,\hat{w}}e_{\hat{w},w'}x_wx_1^{d-2}e_{w',d}\in A,$$
so that $y_wx_wx_1^{d-2}e_{1,d}\in A$ for each $w$ having $2\leq w
\leq d$.  This, together with the previously established
$y_1x_1x_1^{d-2}e_{1,d}\in A$, gives
$$\sum_{w=1}^n y_wx_wx_1^{d-2}e_{1,d} = (\sum_{w=1}^n y_wx_w)x_1^{d-2}e_{1,d}
= 1_K\cdot x_1^{d-2}e_{1,d} \in A$$
so that
$$x_1^{d-2}e_{1,d}\in A.$$
By a procedure analogous to the one we have just completed, which
shows how to obtain $x_1^{d-2}e_{1,d}\in A$ starting from
$x_1^{d-1}e_{1,d}\in A$, we can show that each of the elements
$$x_1^{d-3}e_{1,d},\hspace{.05in} x_1^{d-4}e_{1,d},\hspace{.05in}... \hspace{.05in}, \hspace{.05in}
x_1e_{1,d} \in A,$$
the last of which similarly gives $(\sum_{w=1}^n y_wx_w)e_{1,d} \in
A$, which then finally yields
$$e_{1,d} \in A$$
as desired.   That $e_{d,1}\in A$ follows from Lemma \ref{dual}.
\end{proof}

We finally are in a position to prove the main result of this article.

\begin{theorem}\label{Main}
Let $d,n$ be positive integers, and $K$ any field.  Let $L_{K,n} = L_n$ denote the Leavitt algebra of type $(1,n-1)$
with coefficients in $K$. Then $L_n\cong {\rm M}_d(L_n)$ if and only if ${\rm gcd}(d,n-1)=1$.
\end{theorem}
\begin{proof}  By Proposition \ref{moduletypeofmatrices}, if
${\rm gcd}(d,n-1)>1$ then the module type of ${\rm M}_d(L_n)$ is not $(1,n-1)$, so that ${\rm M}_d(L_n)$ and $L_n$
cannot be isomorphic in this case.

For the implication of interest, suppose ${\rm gcd}(d,n-1)=1$, and suppose $d<n$. By Corollary \ref{method}, we need
only show that the set $A = \{X_1,...,X_n,Y_1,...,Y_n\}$ satisfies the three indicated properties.  That $X_iY_j =
\delta_{i,j}I$ follows directly by the definition of these matrices and Lemma \ref{XYequalsI}.  The equation
$\sum_{j=1}^n Y_jX_j = I$ follows from Lemmas \ref{YiXi1toqplus1}, \ref{IminusEsinA}, and \ref{Listhasdagger}.

For the final property, we must show that $A=<\{X_1,...,X_n,Y_1,...,Y_n\}>={\rm M}_d(L_n)$. It suffices to show that
$x_we_{i,j}\in A$ for all $1\leq w \leq n$ and all $i,j\in \{1,2,...,d\}$, since by Lemma \ref{dual} this will yield
$y_we_{i,j}\in A$ for all $1\leq w \leq n$ and all $i,j\in \{1,2,...,d\}$, and these two collections together clearly
generate all of ${\rm M}_d(L_n)$

By Proposition \ref{possibleplacement} we may assume that the elements from The List have been placed appropriately in
the matrices $X_{q+2},...,X_n$.  Now let $i,j\in \{1,2,...,d\}$. If $i\sim j$ then $e_{i,j}\in A$ by Proposition
\ref{partitionmatrixunits}. So suppose $i\in \hat{S_1}$ and $j\in \hat{S_2}$.  Then $i\sim 1$ and $j\sim d$, so
$e_{i,1}$ and $e_{d,j}$ are each in $A$, again by Proposition \ref{partitionmatrixunits}.  But Proposition \ref{e1dinA}
yields $e_{1,d}\in A$, so that
$$e_{i,1}e_{1,d}e_{d,j}=e_{i,j}\in A.$$
The situation where $i\in \hat{S_2}$ and $j\in \hat{S_1}$ is identical, and thus yields $e_{i,j}\in A$ for all $i,j\in
\{1,2,...,d\}$.   Finally, since each of the elements $\{x_w|1\leq w \leq n\}$ is contained as an entry in one of the
matrices $X_1,...,X_{q+1}$, we can indeed generate all elements of the desired form in $A$.  Thus we have shown that
for ${\rm gcd}(d,n-1)=1$ and $d<n$ we have $L_n\cong {\rm M}_d(L_n)$.

To finish the proof of our main result we need only show that the desired isomorphism holds in case $d\geq n$. Write
$d=q'(n-1)+d'$ with $1\leq d' \leq n-1$.  Then easily ${\rm gcd}(d',n-1)=1$, so the previous paragraph yields $L_n\cong
{\rm M}_{d'}(L_n)$.  But then also $d\equiv d'$ (mod $n-1$), so by Proposition \ref{Basicprops}(1) we get ${\rm
M}_d(L_n)\cong {\rm M}_{d'}(L_n)\cong L_n$, and we are done.
\end{proof}

Notice that Theorem \ref{Main} does not depend on the choice of the positions of the elements from The List in the
non-specified entries of the matrices $X_{q+2}, \dots ,X_n$, other than that the positions are consistent with the
condition allowed by Proposition \ref{possibleplacement}.

\begin{example}
{\rm  We indicated in Section \ref{findmatrices} that $L_5\cong {\rm M}_3(L_5)$; in fact, we provided there five
different sets of appropriate generating matrices of ${\rm M}_3(L_5)$.  Here is yet another set, built by using the
recipe provided in Theorem \ref{Main}.  In this case we have $n=5, d=3, r=2, r-1=1, s=3-1=2, \hat{S_1}=\{1\},
\hat{S_2}=\{2,3\}, S_1=\{1,4\}, S_2=\{2,3,5\}$. The List consists of the $(n-1)(d-1)+1 = 9$ elements
$\{x_1^2,x_2x_1,x_3x_1,x_4x_1,x_5x_1,x_2,x_3,x_4,x_5\}$.  The point to be made here is that only the elements
$x_1^2,x_4x_1$, and $x_4$ can be placed in row $1$ of column $3$, since $S_1=\{1,4\}$. }

$$X_1=\begin{pmatrix}x_{1}&0&0\\
    x_2&0&0\\
    x_3&0&0\end{pmatrix} \hspace{.25in}
 X_2=\begin{pmatrix}x_4&0&0\\
    x_5&0&0\\
    0&1&0\end{pmatrix} \hspace{.25in}$$

 $$X_3=\begin{pmatrix}0&0&x_4\\
    0&0&x_3x_1\\
    0&0&x_2\end{pmatrix} \hspace{.25in}
 X_4= \begin{pmatrix}0&0&x_1^2\\
    0&0&x_2x_1\\
    0&0&x_3\end{pmatrix} \hspace{.25in}
X_5= \begin{pmatrix}0&0&x_4x_1\\
    0&0&x_5\\
    0&0&x_5x_1\end{pmatrix} \hspace{.25in} $$

\end{example}

We finish this section by describing some automorphisms of $L_n$ which arise as a consequence of Theorem \ref{Main}.
There are many possible assignments of the elements on The List to the ``boxes" of the matrices $X_{q+2},X_{q+3}, ...,
X_n$ consistent with the method described in Proposition \ref{possibleplacement}.  In particular, this freedom of
assignment affords an action of the bisymmetric group $S_{d_1}\times S_{d_2}$ on each of the matrices $X_{q+3}, ...,
X_n$ by permuting the entries inside $\hat{S_1}$ and $\hat{S_2}$.   Similarly, we have an action of $S_{e_1}\times
S_{e_2}$ on $X_{q+2}$.  This freedom of assignment also allows an action on each of the $d$ rows in the generating
matrices.  Specifically, for each row $i$ ($1\leq i \leq d$), we can permute the $(i,d)$-entries of the $n-(q+2)$
matrices $X_{q+3}, ..., X_n$; each of the $d\cdot(n-(q+2))!$ such permutations will yield a different set of generators
for ${\rm M}_d(L_n)$. Thus we have described
$$d\cdot (n-(q+2))!e_1!e_2!(d_1!d_2!)^{n-(q+2)}$$
permutations on the entries of the matrices $X_{q+2},X_{q+3}, ..., X_n$, each of which induces a distinct automorphism
of ${\rm M}_d(L_n)$.  In turn, by Theorem \ref{Main}, each then induces an automorphism of $L_n$ whenever ${\rm
gcd}(d,n-1)=1$. These permutations yield automorphisms on $L_n$ which generalize the specific automorphisms of ${\rm
M}_3(L_5)$ described in Section \ref{findmatrices}.

Intriguingly, the types of automorphisms described here and in Section \ref{findmatrices} still do not in general
completely describe all the automorphisms of $L_n$ which arise from producing appropriate sets of generators in ${\rm
M}_d(L_n)$.   We present here two additional specific examples of generating sets inside various-sized matrix rings. In
both cases, the entries used to build the generating matrices $X_1,...,X_n,Y_1,...,Y_n$ are monomials of degree at most
2. In contrast to the previously presented examples, because The List contains monomials of degree up to and including
$d-1$, the examples given here cannot be realized as arising from automorphisms induced by permutations of the entries
of a specific set of generators as constructed in Theorem \ref{Main}.

\medskip

{\bf Example}:  A set of generators of ${\rm M}_4(L_6)\cong L_6$.  (So $d=4$, $d-1=3;$ note that there are no monomials
of degree 3 used in this set.)
$$X_1=\begin{pmatrix}x_{1}&0&0&0\\
    x_2&0&0&0\\
    x_3&0&0&0\\
    x_4&0&0&0\end{pmatrix} \hspace{.25in}
 X_2=\begin{pmatrix}x_5&0&0&0\\
    x_6&0&0&0\\
    0&1&0&0\\
    0&0&1&0\end{pmatrix} \hspace{.25in}
 X_3=\begin{pmatrix}0&0&0&x_1^2\\
    0&0&0&x_2x_1\\
    0&0&0&x_3x_1\\
    0&0&0&x_4x_1\end{pmatrix} \hspace{.25in}$$

 $$X_4=\begin{pmatrix}0&0&0&x_5x_1\\
    0&0&0&x_6x_1\\
    0&0&0&x_1x_2\\
    0&0&0&x_2^2\end{pmatrix} \hspace{.25in}
 X_5=\begin{pmatrix}0&0&0&x_3x_2\\
    0&0&0&x_4x_2\\
    0&0&0&x_5x_2\\
    0&0&0&x_6x_2\end{pmatrix} \hspace{.25in}
X_6=\begin{pmatrix}0&0&0&x_5\\
    0&0&0&x_6\\
    0&0&0&x_3\\
    0&0&0&x_4\end{pmatrix} \hspace{.25in}$$

\medskip

{\bf Example}:  A set of generators of ${\rm M}_5(L_9)\cong L_9$.  (So $d=5$, $d-1=4;$ note that there are no monomials
of degree $3$ or $4$ used in this set.  Also note that, unlike the matrices constructed in Theorem \ref{Main}, there
are entries other than $1_K$ in column $2$.)
$$X_1=\begin{pmatrix}x_{1}&0&0&0&0\\
    x_2&0&0&0&0\\
    x_3&0&0&0&0\\
    x_4&0&0&0&0\\
    x_5&0&0&0&0\end{pmatrix} \hspace{.25in}
 X_2=\begin{pmatrix}x_6&0&0&0&0\\
    x_7&0&0&0&0\\
    x_8&0&0&0&0\\
    x_9&0&0&0&0\\
    0&x_9&0&0&0\end{pmatrix} \hspace{.25in}
 X_3=\begin{pmatrix}0&x_1&0&0&0\\
    0&x_2&0&0&0\\
    0&x_3&0&0&0\\
    0&x_4&0&0&0\\
    0&x_5&0&0&0\end{pmatrix} \hspace{.25in}$$

 $$X_4=\begin{pmatrix}0&x_6&0&0&0\\
    0&x_7&0&0&0\\
    0&x_8&0&0&0\\
    0&0&1&0&0\\
    0&0&0&1&0\end{pmatrix} \hspace{.25in}
X_5=\begin{pmatrix}0&0&0&0&x_1^2\\
    0&0&0&0&x_2x_1\\
    0&0&0&0&x_3x_1\\
    0&0&0&0&x_4x_1\\
    0&0&0&0&x_5x_1\end{pmatrix} \hspace{.25in}
 X_6=\begin{pmatrix}0&0&0&0&x_6x_1\\
    0&0&0&0&x_7x_1\\
    0&0&0&0&x_8x_1\\
    0&0&0&0&x_9x_1\\
    0&0&0&0&x_9\end{pmatrix} \hspace{.25in}$$

$$X_7=\begin{pmatrix}0&0&0&0&x_1x_2\\
    0&0&0&0&x_2^2\\
    0&0&0&0&x_3x_2\\
    0&0&0&0&x_4x_2\\
    0&0&0&0&x_5x_2\end{pmatrix} \hspace{.25in}
 X_8=\begin{pmatrix}0&0&0&0&x_6x_2\\
    0&0&0&0&x_7x_2\\
    0&0&0&0&x_8x_2\\
    0&0&0&0&x_8\\
    0&0&0&0&x_9x_2\end{pmatrix} \hspace{.25in}
 X_9=\begin{pmatrix}0&0&0&0&x_6\\
    0&0&0&0&x_7\\
    0&0&0&0&x_3\\
    0&0&0&0&x_4\\
    0&0&0&0&x_5\end{pmatrix} \hspace{.25in}$$

\bigskip

\section{Applications to C$^*$-algebras and questions about $K_0$}\label{applications}

As mentioned in the Introduction, one consequence of our main result is that we are able to directly and explicitly
establish an affirmative answer to the question posed in \cite{PS}, page 8,  regarding isomorphisms between matrix
rings over Cuntz algebras.

\begin{theorem}\label{C*-fact1}
$ {\rm M}_d(\mathcal{O}_{n}) \cong \mathcal{O}_{n}$ if and only if ${\rm gcd}(d,n-1)=1$.
\end{theorem}
\begin{proof} If ${\rm gcd}(d,n-1)\ne 1$, Then $M_d(\mathcal{O}_n)\not\cong
\mathcal{O}_n$ by \cite{PS}, Corollary 2.4.

So suppose conversely that ${\rm gcd}(d,n-1)=1$. Let $\{ s_1, \dots ,s_n\}\subset \mathcal{O}_n$ be the orthogonal isometries
generating $\mathcal{O}_n$.  These satisfy:
\begin{enumerate}
\item[(i)] For every $1\leq i,j\leq n$, $s_i^*s_j=\delta _{i,j}$, and
\item[(ii)] $1=\sum\limits_{i=1}^ns_is_i^*$.
\end{enumerate}
Now consider the complex Leavitt algebra $L_{\mathbb{C},n}$, and notice that by Proposition \ref{Basicprops}(2) there
exists a (unique) $\mathbb{C}$-algebra morphism
$$\varphi :  L_{\mathbb{C},n}  \rightarrow  \mathcal{O}_n $$
given by the extension of the assignment $x_i  \mapsto  s_i^* $ and $y_i  \mapsto  s_i $ for $1\leq i \leq n$.
 Since $L_{\mathbb{C},n}$ is a
simple algebra, $L_{\mathbb{C},n}\cong \varphi (L_{\mathbb{C},n})$.  But $\mathcal{P}_n=\varphi (L_{\mathbb{C},n})$ is
the complex dense $\ast$-subalgebra of $\mathcal{O}_n$ generated by $\{ s_1, \dots ,s_n\}$ (as a complex algebra). Now
consider the morphism
$$\varphi _d: M_d(L_{\mathbb{C},n})\rightarrow
M_d(\mathcal{O}_n)$$ induced by $\varphi$. (In particular, $\varphi_d(e_{i,j}) = e_{i,j}$ for each matrix unit
$e_{i,j}$, $1\leq i,j \leq d$.)  Notice that, if $X_i$ ($1\leq i\leq n$) is any of the matrices defined in Section
\ref{generators} then, by definition of the elements of The List, $\varphi _d(Y_i)=\varphi _d(X_i)^*$ with respect to
the involution $\ast$ of $M_d(\mathcal{O}_n)$. So, by defining $S_i=\varphi _d(Y_i)$, we get $S^*_i=\varphi _d(X_i)$,
and thus $\{ S_1, \dots ,S_n\}\subset M_d(\mathcal{O}_n)$ is a family of $n$ orthogonal isometries satisfying
$I_d=\sum\limits_{i=1}^nS_iS_i^*$. Hence, by \cite{Cu}, Theorem 1.12, there exists an isomorphism
$$\Phi :\mathcal{O}_n\rightarrow C^*(S_1, \dots , S_n)\subseteq M_d(\mathcal{O}_n)$$
defined by the rule $\Phi (s_i)=S_i$ for every $1\leq i\leq n$. Now,
applying Theorem \ref{Main} to $\mathcal{P}_n$ and
$M_d(\mathcal{P}_n)$ (via $\varphi$), for every $1\leq i,j\leq d$
and for every $1\leq k\leq n$ we have
$$s_ke_{i,j}=\varphi_d(y_ke_{i,j})\in C^*(S_1, \dots , S_n),$$
so that the generators of $M_d(\mathcal{O}_n)$ lie in $C^*(S_1, \dots , S_n)$. Thus, $\mathcal{O}_n\cong
M_d(\mathcal{O}_n)$ via $\Phi$, so we are done.
\end{proof}

As mentioned previously, the affirmative answer to the isomorphism question for matrix rings over Cuntz algebras
provided in Theorem \ref{C*-fact1} is indeed already known, a byproduct of \cite{Ph1}, Theorem 4.3(1). However, the
method we have provided in Theorem \ref{C*-fact1} is significantly more elementary, and provides an explicit
description of the germane isomorphisms (such an explicit description has previously not been known).

\smallskip

A second interesting consequence of Theorem \ref{Main} is that the class of matrices over Leavitt algebras is
classifiable using K-theoretic invariants.  (For additional information about purely infinite simple algebras and their
$K$-theory, see \cite{AGP}.)

\begin{theorem}\label{C*-fact2}
Let $\mathcal{L}$ denote the set of purely infinite simple
$K$-algebras
$$\{{\rm M}_d(L_n) | d,n \in \mathbb{N}\}.$$
Let $B,B'\in \mathcal{L}$. Then $B\cong B'$ if and only if there
is an isomorphism $\phi: K_0(B)\rightarrow K_0(B')$ for which
$\phi([1_B])=[1_{B'}]$.
\end{theorem}
\begin{proof}
It is well known (see e.g. \cite {Ros}, page 5) that any unital isomorphism $f: B\rightarrow B'$ induces a group
isomorphism $K_0(f): K_0(B)\rightarrow K_0(B')$ sending $[1_B]$ to $[1_{B'}]$.

To see the converse, first notice that, for any $B\in \mathcal{L}$, $B={\rm M}_d(L_n)$ for suitable $d, n\in \N$. It is
well known that
$$(K_0({\rm M}_d(L_n)), [1_{{\rm M}_d(L_n)}])\cong (\Z/(n-1)\Z, [d])$$
(see e.g \cite{B} or \cite{AGP}). Hence, if $B'={\rm M}_k(L_m)$ for suitable $k, m\in \N$, then the existence of an
isomorphism $\phi: K_0(B)\rightarrow K_0(B')$ forces that $n=m$.

Now, since every automorphism of $\Z/(n-1)\Z$ is given by multiplication by an element $1\leq l\leq n-1$ such that
${\rm gcd}(l,n-1)=1$, the hypothesis $\phi([1_B])=[1_{B'}]$ yields that $[k]=[dl]\in \Z/(n-1)\Z$, i.e., that $k\equiv
dl$ (mod $n-1$). So Proposition \ref{Basicprops}(1) gives that
$${\rm M}_k(L_n)\cong {\rm M}_{dl}(L_n)\cong {\rm M}_d({\rm M}_l(L_n)).$$
Since ${\rm gcd}(l,n-1)=1$, we have ${\rm M}_l(L_n)\cong L_n$ by Theorem \ref{Main}. Hence, ${\rm M}_d({\rm
M}_l(L_n))\cong {\rm M}_d(L_n)$, whence
$${\rm M}_k(L_n)\cong {\rm M}_{dl}(L_n)\cong {\rm M}_d({\rm M}_l(L_n))\cong {\rm M}_d(L_n),$$ as
desired.
\end{proof}

A significantly  more general C$^*$-algebraic analog of Theorem \ref{C*-fact2} is well-known for the class of unital
purely infinite simple C*-algebras, as a consequence of the  powerful work  of Kirchberg and Phillips, \cite{K} and
\cite{P}. However, even in the concrete case of the subclass $\{{\rm M}_d(\mathcal{O}_n) | d,n \in \mathbb{N}\}$, the
existence of the previously known isomorphisms in the C$^*$-algebra setting (to wit, the aforementioned  results of
R{\o}rdam, Kirchberg and Phillips) depend on deep results which produce no explicit isomorphisms. A natural question in
this context is whether \cite{P}, Theorem 4.2.4, has an algebraic counterpart. In \cite{AAP} the authors establish a
partial affirmative answer to this question for a large class of purely infinite simple algebras.

\section{Graded isomorphisms between Leavitt algebras and their matrix rings.}\label{graded}

In this final section we incorporate the natural ${\mathbb Z}$-grading on the Leavitt algebras into our analysis.  As
one consequence, we will show that the sets of matrices which arise in the proof of \cite{L1}, Theorem 5, cannot in
general generate ${\rm M}_d(L_n)$.

The ${\mathbb Z}$-grading on $L_{K,n}$ is given as follows.  We define the degree of a monomial of the form
$y_i^tx_j^u$ by setting
$$deg(y_i^tx_j^u) = u-t,$$
 and extending linearly to all of $L_{K,n}$.  This is precisely the ${\mathbb Z}$-grading on $L_{K,n}$ induced by setting $deg(X_i)=1, deg(Y_i)=-1$
 in $R = K<X_1,...,X_n,Y_1,...,Y_n>$, and then grading the factor ring $L_n = R/I$ in the natural way.
 (We note that the relations which define $L_n$  are homogeneous in this grading of $R$.)

 It was shown in
 \cite{AAn1} that, in this grading, $(L_n)_0 \cong \lim \limits_{\longrightarrow}{}_{t\in \N}({\rm M}_{n^t}(K))$.
 Here the connecting homomorphisms are unital (so that the direct limit is unital); the homomorphism from ${\rm M}_{n^t}(K)$
 to ${\rm M}_{n^{t+1}}(K)$ is given by sending any matrix of the form $(a_{i,j})$ to the matrix $(a_{i,j}I_n)$.

We will need the following easily proved result about unital direct limits of rings.   For a unital ring $R$, we say
that a finite set $E = \{e_1,...,e_p\}$ of idempotents in $R$ is {\it complete, orthogonal, pairwise isomorphic} in
case $1_R = e_1 + ... + e_p$, $e_ie_j=0$ for all $i\neq j$, and  $Re_i\cong Re_j$ as left $R$-modules for all $1\leq
i,j \leq p$. In particular, in this situation we have $R\cong \oplus_{i=1}^p Re_i$ as left $R$-modules.

\begin{lemma}\label{isomorphicinsubrings}
Suppose $R$ is a unital direct limit of rings $R =\lim \limits_{\longrightarrow}{}_{t\in \N}(R_t)$ (so we are assuming
that connecting homomorphism $R_t \rightarrow R_{t+1}$ is unital for each $t\in \N$). Suppose $R$ contains a complete
orthogonal pairwise isomorphic set of $p$ idempotents. Then there exists $m\in \N$ so that $R_m$ contains a complete
orthogonal pairwise isomorphic set of $p$ idempotents.
\end{lemma}
\begin{proof}
Let $E = \{e_1,...,e_p\}$ denote the indicated set in $R$.   It is well known (see e.g. \cite{J}, Proposition III.7.4)
that for idempotents $e$ and $f$ in any ring $R$, $Re\cong Rf$ as left $R$-modules if and only if there exist elements
$x,y$ in $R$ such that $x=exf$, $y=fye$, $xy=e$, and $yx=f$.   For each two-element subset $\{e_i,e_j\}$ of $E$ let
$\{x_{i,j},y_{i,j}\}$ denote a pair of associated elements whose existence is ensured by the supposed isomorphism
$Re_i\cong Re_j$.   Now pick $m\in \N$ with the property that $R_m$ contains the finite set $\{x_{i,j},y_{i,j}\mid
1\leq i,j \leq p\}$; such $m$ exists by definition of direct limit.   Then necessarily $R_m$ contains $E$, as
$x_{i,j}y_{i,j}=e_i$ for each $1\leq i \leq p$.  Now invoking the previously cited result from \cite{J}, and using the
hypothesis that the direct limit has unital connecting homomorphisms, we conclude that $E$ is a complete orthogonal
pairwise isomorphic set of $p$ idempotents in $R_m$.
\end{proof}

\begin{lemma}\label{matrixsizes}
Let $S$ be any unital ring, let $K$ be a field, and let $p$ be any positive integer.
\begin{enumerate}
\item If $p\mid d$, then the matrix ring $T={\rm M}_d(S)$ contains
a complete, orthogonal, pairwise isomorphic set of $p$ idempotents.
\item If the matrix ring $T={\rm M}_d(K)$ contains a complete, orthogonal, pairwise isomorphic set of $p$ idempotents,
then $p\mid d$.
\end{enumerate}
\end{lemma}
\begin{proof}  For (1), writing $d=pq$ and using the isomorphism ${\rm M}_d(S) \cong {\rm M}_p({\rm M}_q(S))$ produces such a
set, where we take $E$ to be the set of $p$ matrix idempotents in ${\rm M}_p({\rm M}_q(S))$.

For (2), let $E$ be such a set.  The ring $T={\rm M}_d(K)$ is semisimple artinian, with composition length $d$.  As the
left $T$-modules $Te_i$ generated by the elements of $E$ are pairwise isomorphic, each must have the same composition
length, which we denote by $q$. But $T\cong \oplus_{i=1}^p Te_i$, which yields that $pq=d$.
\end{proof}

With these two lemmas in hand, we are ready to prove the main result of this section.

\begin{proposition}\label{gradediso}
The algebras  $L_n$ and ${\rm M}_d(L_n)$ are isomorphic as $\mathbb{Z}$-graded algebras if and only if there exists
$\alpha \in \mathbb{N}$ such that $d \mid n^{\alpha}$.
\end{proposition}
\begin{proof}
First suppose there exists $\alpha \in \mathbb{N}$ such that $d \mid n^{\alpha}$.  Then the explicit isomorphism
provided in \cite{PS}, Proposition 2.5,  between the indicated matrix rings over Cuntz algebras is easily seen to
restrict to an isomorphism of the analogously-sized matrix rings over Leavitt algebras.  Furthermore, the isomorphism
preserves the appropriate grading on these algebras, thus yielding the first implication.  (For clarity, an explicit
example of this isomorphism in a particular case is given below.)

Conversely, suppose the algebras  $L_n$ and ${\rm M}_d(L_n)$ are isomorphic as $\mathbb{Z}$-graded algebras.  Then
necessarily the $0$-components of these algebras are isomorphic.  It is easy to show that the $0$-component of ${\rm
M}_d(L_n)$ is isomorphic to ${\rm M}_d (\lim \limits_{\longrightarrow}{}_{t\in \N}(M_{n^t}(K)))$. Now let $p$ be any
prime number with $p\mid d$.  Then by Lemma \ref{matrixsizes}(1), ${\rm M}_d (\lim \limits_{\longrightarrow}{}_{t\in
\N}(M_{n^t}(K)))$ contains a complete orthogonal pairwise isomorphic set of $p$ idempotents.  Using the isomorphism
between $0$-components, we get a complete orthogonal pairwise isomorphic set of $p$ idempotents in $\lim
\limits_{\longrightarrow}{}_{t\in \N}({\rm M}_{n^t}(K))\cong (L_n)_0$. But by Lemma \ref{isomorphicinsubrings}, this
implies that there exists an integer $u$ so that the matrix ring ${\rm M}_{n^u}(K)$ contains a complete orthogonal
pairwise isomorphic set of $p$ idempotents. By Lemma \ref{matrixsizes}(2) this implies that $p\mid n^u$, so that $p\mid
n$ as $p$ is prime.    Thus we have shown that any prime $p$ which divides $d$ also necessarily divides $n$, so that
$d$ indeed divides some power of $n$ as desired.
\end{proof}

\begin{corollary}\label{gradedisoonlywhenddividesnalpha}
Suppose ${\rm gcd}(d,n-1)=1$.  Suppose $W = \{X_1,...,X_n,Y_1,...,Y_n\}$ is a set of $2n$ matrices in ${\rm M}_d(L_n)$
which satisfy the conditions of Proposition \ref{Basicprops}(2).  Suppose further that each entry of $X_i$ (resp.
$Y_i$)  is either $0$ or a monomial of degree 1 (resp. degree -1).   If $W$ generates ${\rm M}_d(L_n)$ as a
$K$-algebra, then
 $d\mid n^{\alpha}$ for some positive integer $\alpha$.

 In particular, let $W$ be the set of $2n$ matrices
$\{X_1,...,X_n,Y_1,...,Y_n\}$ constructed in \cite{L1}, Theorem 5. Then $W$ generates ${\rm M}_d(L_n)$ as a $K$-algebra
if and only if $d\mid n^{\alpha}$ for some positive integer $\alpha$.
 \end{corollary}
\begin{proof}
If $W$ satisfies the indicated conditions, then the homomorphism from $L_n$ to ${\rm M}_d(L_n)$ induced by the
assignment $x_i \mapsto X_i$ and $y_i \mapsto Y_i$ in fact would be a graded isomorphism, and the result follows from
Proposition \ref{gradediso}.

In the specific case of the $2n$ matrices described in \cite{L1}, Theorem 5,  the matrices are of the indicated type,
and were shown in \cite{PS} to generate ${\rm M}_d(L_n)$.
 \end{proof}

It is instructive to compare and contrast the two types of generating sets of ${\rm M}_d(L_n)$ which can be constructed
in case $d\mid n^{\alpha}$ for some $\alpha$. Let $d=3,n=6$. Here are the six matrices $\{X_1,...,X_6\}$ which arise in
the aforementioned construction presented in \cite{PS}.

$$X_1=\begin{pmatrix}x_{1}&0&0\\
    x_2&0&0\\
    x_3&0&0\end{pmatrix} \hspace{.25in}
 X_2=\begin{pmatrix}x_4&0&0\\
    x_5&0&0\\
    x_6&0&0\end{pmatrix} \hspace{.25in}
 X_3=\begin{pmatrix}0&x_1&0\\
    0&x_2&0\\
    0&x_3&0\end{pmatrix} \hspace{.25in}$$

 $$X_4=\begin{pmatrix}0&x_4&0\\
    0&x_5&0\\
    0&x_6&0\end{pmatrix} \hspace{.25in}
 X_5= \begin{pmatrix}0&0&x_1\\
    0&0&x_2\\
    0&0&x_3\end{pmatrix} \hspace{.25in}
X_6= \begin{pmatrix}0&0&x_4\\
    0&0&x_5\\
    0&0&x_6\end{pmatrix} \hspace{.25in} $$

In particular, all of these are of degree $1$ in the ${\mathbb Z}$-grading, so that the assignment $x_i \mapsto X_i$
(and $x^*_i \mapsto X^*_i$) from $L_6$ to ${\rm M}_3(L_6)$ extends to a graded homomorphism, which can be shown in a
straightforward way (using the argument given in \cite{PS}) to be a graded isomorphism.

\smallskip

In contrast, we now present one (of many) sets of generators for ${\rm M}_3(L_6)$ which arises from our construction.
When $n=6,d=3$ then the appropriate data from our main result are as follows: $6 = 1\cdot 3 + 3$, so $r=3$, $r-1 = 2$,
$s=3-2 = 1$, $\hat{S_1}=\{1,2\}$, $\hat{S_2} = \{3\}$, $S_1=\{1,2,4,5\}$, $S_2 = \{3,6\}$.   So one possible collection
of appropriate generating matrices in ${\rm M}_3(L_6)$ is

$$X_1=\begin{pmatrix}x_{1}&0&0\\
    x_2&0&0\\
    x_3&0&0\end{pmatrix} \hspace{.25in}
 X_2=\begin{pmatrix}x_4&0&0\\
    x_5&0&0\\
    x_6&0&0\end{pmatrix} \hspace{.25in}
 X_3=\begin{pmatrix}0&1&0\\
    0&0&x_1^2\\
    0&0&x_3x_1\end{pmatrix} \hspace{.25in}$$

 $$X_4=\begin{pmatrix}0&0&x_2x_1\\
    0&0&x_4x_1\\
    0&0&x_6x_1\end{pmatrix} \hspace{.25in}
 X_5= \begin{pmatrix}0&0&x_5x_1\\
    0&0&x_2\\
    0&0&x_3\end{pmatrix} \hspace{.25in}
X_6= \begin{pmatrix}0&0&x_4\\
    0&0&x_5\\
    0&0&x_6\end{pmatrix} \hspace{.25in} $$

\medskip

We close this article by providing a brief historical perspective on this question.   As mentioned earlier, Leavitt
showed in \cite{L1} that if $R$ has module type $(1,n-1)$, then ${\rm M}_d(R)$ has module type $(1, \frac{n-1}{{\rm
gcd}(d,n-1)}$). The validity of this result is justified by the presentation of an appropriate set of elements inside
${\rm M}_d(R)$. In the situation where $R = L_n$ and ${\rm gcd}(d,n-1)=1$, it turns out that the appropriate set of
elements inside ${\rm M}_d(R)$ is simply a lexicographic ordering of the variables $\{x_1,...,x_n,y_1,...,y_n\}$, using
a straightforward algorithm.  (An example of this process was given in Section \ref{findmatrices}.) In the particular
case when $d\mid n^{\alpha}$, the set of elements so constructed coincides with the set of elements analyzed by Paschke
and Salinas in \cite{PS}; furthermore, this set just happens to generate all of ${\rm M}_d(L_n)$. However, as noted in
Corollary \ref{gradedisoonlywhenddividesnalpha}, the analogous set of elements cannot generate all of ${\rm M}_d(L_n)$
when $d$ is not a divisor of some power of $n$.   Thus, in order to establish our main result (Theorem \ref{Main}), it
was necessary to build a completely different set of tools than those which had already been used in this arena.

Corollary \ref{gradedisoonlywhenddividesnalpha} shows that in general we cannot find generating sets of size $2n$
inside ${\rm M}_d(L_n)$ in which each of the entries in the $n$ matrices has degree $1$ (resp., each of the entries in
the $n$ dual matrices has degree -1).   In our main result we have shown that we can find generating sets of size $2n$
inside ${\rm M}_d(L_n)$ in which each of the entries in the $n$ matrices has degree less than or equal to $d-1$ (resp.,
each of the entries in the $n$ dual matrices has degree greater than or equal to $1-d$).  Reflecting on the examples
given at the end of Section \ref{maintheorem}, it would be interesting to know whether in general it is possible to
find generating sets of size $2n$ inside ${\rm M}_d(L_n)$ in which each of the entries in the $n$ matrices has degree
less than or equal to $2$ (resp. each of the entries in the $n$ dual matrices has degree greater than or equal to
$-2$).

\bibliographystyle{amsplain}

\end{document}